\documentclass[11pt,fleqn,a4paper]{article}
\usepackage{amssymb,latexsym,amsmath,amsfonts}
\usepackage{verbatim}
\usepackage{graphicx}

    \numberwithin{equation}{section}

    \def\Re{{\rm Re \,}}
    \def\Im{{\rm Im \,}}

\DeclareMathOperator{\Ai}{Ai}
\DeclareMathOperator{\diag}{diag}

\DeclareMathOperator{\supp}{supp}

    \newtheorem{theorem}{Theorem}[section]
    \newtheorem{lemma}[theorem]{Lemma}
    
    \newtheorem{proposition}[theorem]{Proposition}
    
    \newtheorem{Definition}[theorem]{Definition}
    
    \newtheorem{Remark}[theorem]{Remark}
    
    \newtheorem{Example}[theorem]{Example}

    \newenvironment{proof}%
    {\rm \trivlist \item[\hskip \labelsep{\bf Proof. }]}%
    {\hspace*{\fill}$\Box$\endtrivlist}
    {\rm \trivlist \item[\hskip \labelsep{\bf Proof}]}%
    {\hspace*{\fill}$\Box$\endtrivlist}

\begin{document}
 \begin{center} \Large\bf
    Asymptotics of non-intersecting Brownian motions
    and a $4 \times 4$ Riemann-Hilbert
    problem
\end{center}

    \begin{center}  \large
        E. Daems, A.B.J. Kuijlaars, and W. Veys
        \footnote{
        The first two authors are supported by FWO-Flanders project G.0455.04,
by K.U. Leuven research grant OT/04/24, and by the
 European Science Foundation Program MISGAM. The second author is
 also
supported by INTAS Research Network 03-51-6637,
and by a grant from the Ministry of Education and
Science of Spain, project code MTM2005-08648-C02-01.
The third author is supported by FWO-Flanders project G.0318.06.}
 \\[3ex]
           \normalsize \em
            Department of Mathematics, Katholieke Universiteit Leuven, \\
            Celestijnenlaan 200 B, 3001 Leuven, Belgium \\[1ex]
            \rm evi.daems@wis.kuleuven.be \\
            \rm arno.kuijlaars@wis.kuleuven.be \\
            \rm wim.veys@wis.kuleuven.be \\[1ex]
     \end{center}\ \\[1ex]

\begin{abstract}
We consider $n$ one-dimensional Brownian motions, such that $n/2$
Brownian motions start at time $t=0$ in the starting point $a$ and
end at time $t=1$ in the endpoint $b$ and the other $n/2$ Brownian
motions start at time $t=0$ at the point $-a$ and end at time
$t=1$ in the point $-b$, conditioned that the $n$ Brownian paths
do not intersect in the whole time interval $(0,1)$. The
correlation functions of the positions of the non-intersecting
Brownian motions have a determinantal form with a kernel that is expressed
in terms of multiple Hermite polynomials of mixed type. We
analyze this kernel in the large $n$ limit
for the case $ab<1/2$. We find  that the limiting mean density
of the positions of the Brownian motions is supported on one or two intervals
and that the correlation kernel has the usual scaling limits from random
matrix theory, namely the sine kernel in the bulk and the Airy kernel
near the edges.
\end{abstract}

\bigskip

\bigskip

{\em Keywords:}
Multiple orthogonal polynomials, non-intersecting Brownian motion,
Riemann-Hilbert problem, steepest descent analysis

\newpage

\section{Introduction}
\label{section1}

\begin{figure}
    \begin{center}
        \includegraphics[height=4.5cm]{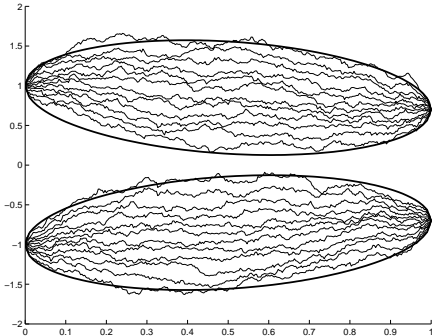}
        \caption{\label{figuur1b}
        Non-intersecting Brownian motions when $ab>1/2$. Here we have
        chosen $a=1$ and $b=0.7$. For $ab > 1/2$
        we have that the two groups of Brownian motions remain separated during
        the full time interval $(0,1)$. }
    \end{center}
\end{figure}

\begin{figure}
    \begin{center}
        \includegraphics[height=4.5cm]{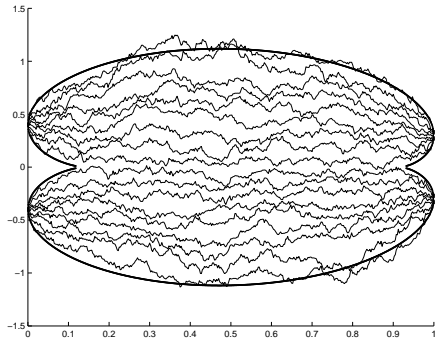}
        \caption{ \label{figuur1a}
        Non-intersecting Brownian motions when $ab <1/2$. Here we have
        chosen $a=0.4$ and $b=0.3$. For $ab < 1/2$ we have that the two groups of
        Brownian paths starting at $\pm a$ come together at the first critical
        time $t_{c,1}$, they merge and continue as one group until the second
        critical time $t_{c,2}$, after which they split again and end at $\pm b$.}
    \end{center}
\end{figure}

\begin{figure}
    \begin{center}
        \includegraphics[height=4.2cm]{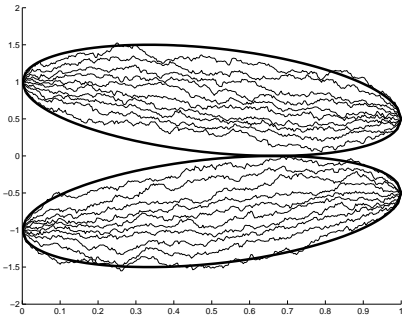}
        \caption{\label{figuur1c}
        Non-intersecting Brownian motions when $ab=1/2$. Here we have
        chosen $a=1$ and $b=1/2$. For $ab=1/2$, we have that the two groups of
        Brownian paths touch each other at a critical time $t_c$.}
    \end{center}
\end{figure}

Consider $n$ one-dimensional Brownian motions conditioned
not to intersect in the time interval $(0,1)$. We assume that
$n$ is even and that $n/2$ Brownian motions start at the
position $a>0$ at time $t=0$ and end at $b>0$ at time $t=1$,
while $n/2$ Brownian motions start at $-a$ and end at $-b$.
If we let $n \to \infty$ and at the same time rescale the
variance of the Brownian motion with a factor $1/n$, then the
Brownian motions fill out a region in the $tx$-plane that
looks like one of the regions shown in Figures \ref{figuur1b}--\ref{figuur1c}.

The shape of the region depends on the product $ab$. For $ab$
greater than some critical value (which in the units that we
will be using is $1/2$) the starting positions
$\pm a$ and the end positions $\pm b$ are sufficient apart
so that the two groups of paths remain separated. The Brownian
motions then fill out two ellipses as can be seen in Figure~\ref{figuur1b}.
Each group essentially behaves like $n/2$ non-intersecting
Brownian motions with a single starting and
end position. This is a variation of Dyson's Brownian motion for
the behavior of the eigenvalues of a Hermitian matrix whose elements
evolve according to a Brownian motion \cite{Dyson}.
In that model, it holds that at each time $t \in (0,1)$ the positions
of the paths are distributed like the eigenvalues of a GUE matrix and
so follow Wigner's semi-circle law as $n \to \infty$, see \cite{Mehta}.
For the present model and $ab > 1/2$, we then also expect
to find at each time $t \in (0,1)$ the scaling limits
for the correlation functions that are known from
random matrix theory, and that are expressed in terms of
the sine kernel in the bulk, and the Airy kernel at the edge points.

The behavior is different for $ab$ less than the critical value $1/2$.
In that case the starting and end positions are not that far apart
and the two groups of paths will interact with each other.
There are two critical times $t_{c,1}$ and $t_{c,2}$ such that the
two groups are separated up to the first critical time $t_{c,1}$.
At $t_{c,1}$ they  merge and continue as
one group of paths until the second critical time $t_{c,2}$ when they
split again, see Figure~\ref{figuur1a}. This is an extension
of the case of two starting positions and one end position
which was studied in \cite{ABK,BK2,BK3} with the use of Riemann-Hilbert
techniques and in \cite{TW1} by classical steepest descent techniques.
There it was found that
for each time $t \in (0,1)$, the scaling limits are still expressed
in terms of the sine kernel in the bulk and in terms of the Airy kernel
at the edge, except at the cusp point at the critical time. At
the cusp the scaling limits are expressed in terms
of Pearcey kernels \cite{BK3,TW1} that were first described by
Br\'ezin and Hikami \cite{BH1,BH2}. Tracy and Widom \cite{TW1}
also identified a Pearcey process that is further
discussed in \cite{AvM1,AvM2,OR} as well.

Returning to the model with two starting positions and two end positions
we have a critical separation $ab = 1/2$. Then the starting
and end positions have a critical separation so that the two
groups of paths just touch in one point, see Figure~\ref{figuur1c}.
Here we expect new critical behavior that can be expressed in terms
of an as yet unknown kernel.

It is the aim of the present paper to treat the case $ab < 1/2$,
$t \neq t_{c1},t_{c,2}$ with the methods of \cite{ABK,BK2,BK3}.
That is, we use a steepest descent analysis of a relevant
Riemann-Hilbert problem that was given by the first two authors in \cite{DK2}.
The Riemann-Hilbert problem has size $4 \times 4$ and its solution is
constructed out of multiple Hermite polynomials of mixed type.

As a result of the asymptotic analysis we find the usual
sine kernel in the bulk and the Airy kernel at the edge points,
see Theorems \ref{theorem22} and \ref{theorem23} below,
thereby providing further evidence for the universality of these
kernels. We are not aware of a double integral representation of
the finite $n$ correlation kernels, so that an asymptotic analysis
of integrals as in \cite{TW1} may not be possible in this case.

The connections between random matrices and non-intersecting
Brownian motions are well-known, see for instance \cite{Joh}
and \cite{TW2} for a very recent paper on non-intersecting
Brownian excursions.
Unfortunately we do not know if there exists a corresponding
random matrix model for the case studied in this paper.

\section{Statement of results} \label{section2}

\subsection{Correlation kernel}
It follows from the classical paper of Karlin and McGregor \cite{KMcG}
that the (random) positions of the Brownian motions at time $t \in (0,1)$
are a determinantal point process. This means that there
is a kernel $K_n$ so that for each $m$
we have that the $m$-point correlation function
\[ R_m(x_1,\ldots, x_m) = \frac{n!}{(n-m)!}
    \int \ldots \int p_{n,t}(x_1, \ldots, x_n) dx_{m+1} \ldots dx_n, \]
where $p_{n,t}(x_1,\ldots, x_n)$ denotes the joint probability density
function for the positions of the paths at time $t$, is given
by the determinant
\[ R_m(x_1, \ldots, x_m) = \det(K_n(x_i,x_j))_{i,j=1,\ldots,m}, \]
see \cite{Sosh}.

Indeed, we have by \cite{KMcG} that
\begin{equation} \label{KarlinMcGregor}
    p_{n,t}(x_1, \ldots, x_n) \propto \det(P_n(t,a_j,x_k))_{j,k=1,\ldots,n}
    \det(P_n(1-t,x_k,b_j))_{j,k=1,\ldots,n}
\end{equation}
if the non-intersecting paths start at $a_1 < a_2 < \cdots < a_n$ at time $t=0$
and end at $b_1 < b_2 < \cdots < b_n$ at time $t=1$, where $P_n(t,a,x)$
is the transition probability for Brownian motion with variance $1/n$,
\[ P_n(t,a,x) = \frac{\sqrt{n}}{\sqrt{2\pi t}} e^{-\frac{n}{2t}(x-a)^2}. \]
Thus (\ref{KarlinMcGregor}) is a biorthogonal ensemble \cite{Bor},
which is a special case of a determinantal point process. These
ensembles are also further studied in connection with random
matrix theory in  \cite{DF}. The correlation kernel is
obtained by biorthogonalizing the two sets of functions
$P_n(t,a_j,\cdot)$ and $P_n(1-t,\cdot,b_j)$ which results in
functions $\phi_j$ and $\psi_j$ for $j=1,\ldots, n$, say, and then
putting
\begin{equation} \label{Knseries}
    K_n(x,y) = \sum_{j=1}^n \phi_j(x) \psi_j(y).
\end{equation}

In the situation of the present paper we have to take the confluent limit
$a_j \to a$, $b_j \to b$ for $j=1,\ldots, n/2$, and $a_j \to -a$, $b_j \to -b$
for $j=n/2+1, \ldots, n$. Then we continue to have a biorthogonal ensemble
and the structure of the kernel (\ref{Knseries}) remains the same.
In \cite{DK2} the functions $\phi_j$ and $\psi_j$ are written in
terms of certain polynomials that were called multiple Hermite
polynomials of mixed type. We will use here only the fact that
$K_n$ is expressed in terms of the solution of a Riemann-Hilbert
problem, see \cite{DK2} and (\ref{defKinY}) below.

We assume throughout that
\begin{align} \label{assume}
    0 < ab < 1/2 \qquad \mbox{ and } \qquad
    t \in (0,1) \setminus \{ t_{c,1}, t_{c,2}\}.
    \end{align}
So we exclude the critical times from our considerations.

\subsection{Limiting mean density}
Our first result deals with the limiting mean density
    \begin{align} \label{rhox}
        \rho(x) = \lim_{n \to \infty} \frac{1}{n}K_n(x,x).
    \end{align}
\begin{theorem}\label{theorem21}
    Let $a, b > 0$ and $t \in (0,1)$ such that {\rm (\ref{assume})} holds.
    Then the limiting mean density {\rm (\ref{rhox})} of the positions of
    the Brownian paths at time $t$ exists. It is supported
    on one interval $[-z_1,z_1]$ if $t \in (t_{c,1}, t_{c,2})$
    and on two intervals $[-z_1,-z_2] \cup [z_2,z_1]$ if $t \in (0,t_{c,1}) \cup (t_{c,2},1)$.
    In all cases the density $\rho(x)$ is expressed as
   \begin{align}
       \rho(x) = \frac{1}{\pi}|\Im \xi(x)|,
   \end{align}
   where $\xi=\xi(x)$ is a solution of the    equation
\begin{align}\label{algemeenriemannoppervlak}
    \begin{split}
     &\xi^4-\frac{2z}{t(1-t)}\xi^3+\left(\frac{z^2}{t^2(1-t)^2}-\frac{a^2}{t^2}-\frac{b^2}{(1-t)^2}
    + \frac{1}{t(1-t)}\right)\xi^2 \\ &+ \left(\frac{2b^2}{t(1-t)^3}-\frac{1}{t^2(1-t)^2}\right)z\xi
    -\frac{b^2z^2}{t^2(1-t)^4}=0.
    \end{split}
\end{align}
    The density $\rho$ is real and analytic on the interior of its support
    and it vanishes like a square root at the edges of its support, i.e.,
    there exists a constant $c_1$ such that
    \begin{align}\label{rhoz1}
        \rho(x) = \frac{c_1}{\pi}|x\mp z_1|^{1/2}\left(1+o(1)\right) \quad \mbox{ as
        $x \to \pm z_1$, $x \in \supp \rho$},
    \end{align}
    and, in case $t \in (0,t_{c,1}) \cup (t_{c,2},1)$, there exists
    a constant $c_2$ such that
    \begin{align} \label{rhoz2}
        \rho(x) = \frac{c_2}{\pi}|x \mp z_2|^{1/2}(1+o(1))
        \quad \mbox{ as $x \to \pm z_2$, $x \in
        \supp \rho$}.
    \end{align}
\end{theorem}

\subsection{Scaling limits of the kernel}
As in \cite{ABK} and \cite{BK2}, the local eigenvalue results are formulated in
terms of a rescaled version of the kernel $K_n$,
\begin{align}\label{rescaledkernel2}
    \hat{K}_n(x,y) = e^{n(h(x)-h(y))} K_n(x,y)
\end{align}
for some function $h$. Note that this change in the kernel does not
affect the determinants
\[ \det(K_n(x_i,x_j))_{1\leq i, j \leq m} \]
that give the correlation functions for the determinantal process.

Theorems \ref{theorem22} and \ref{theorem23} show that the
kernel has the scaling limits that are universal for unitary
ensembles in random matrix theory \cite{DKMVZ1},
namely the sine kernel in the bulk, and the Airy kernel at the edge.
\begin{theorem}\label{theorem22}
    Let $a, b > 0$ and $t \in (0,1)$ such that {\rm (\ref{assume})} holds.
    Let $z_1, z_2$ be as in Theorem {\rm \ref{theorem21}}.
    Then there
    exists a function $h$ such that the following
    holds for the rescaled kernel {\rm(\ref{rescaledkernel2})}.
    For every  $x_0 \in (-z_1,z_1)$ (in case $t \in (t_{c,1}, t_{c,2})$)
    or $x_0 \in (-z_1,-z_2) \cup (z_2, z_1)$ (in case $t \in (0,t_{c,1}) \cup (t_{c,2},1)$)
    and $u,v \in \mathbb R$, we have
        \begin{align} \label{sinekernellimit}
            \lim_{n \to
            \infty}\frac{1}{n\rho(x_0)}\hat{K}_n\left(x_0+\frac{u}{n\rho(x_0)},
            x_0+\frac{v}{n\rho(x_0)}\right) = \frac{\sin
            \pi(u-v)}{\pi(u-v)}.
        \end{align}
\end{theorem}

\begin{theorem}\label{theorem23}
    Let $a, b > 0$ and $t \in (0,1)$ such that {\rm (\ref{assume})} holds.
    Let $z_1, z_2$ be as in Theorem {\rm \ref{theorem21}}.
    Then there exists a function $h$ such that the following
    holds for the rescaled kernel {\rm(\ref{rescaledkernel2})}.
    For every $u, v\in \mathbb R$, we have
        \begin{align} \nonumber  & \lim_{n \to \infty}
            \frac{1}{(c_1n)^{2/3}}
            \hat{K}_n\left( z_1 + \frac{u}{(c_1n)^{2/3}},
             z_1 +\frac{v}{(c_1n)^{2/3}}\right) \\ &\qquad \qquad =
            \frac{\Ai(u)\Ai'(v)-\Ai'(u)\Ai(v)}{u-v},
            \label{airykernellimit1}
        \end{align}
    where $\Ai$ is the usual Airy function, and $c_1$
    is the constant defined in {\rm(\ref{rhoz1})}.

    In addition, if $t \in (0,t_{c,1}) \cup (t_{c,2},1)$,
    then for every $u, v \in \mathbb R$, we have
        \begin{align} \nonumber
            & \lim_{n \to \infty}
            \frac{1}{(c_2n)^{2/3}}
            \hat{K}_n\left( z_2 - \frac{u}{(c_2n)^{2/3}},
             z_2 -\frac{v}{(c_2n)^{2/3}}\right) \\&\qquad \qquad =
            \frac{\Ai(u)\Ai'(v)-\Ai'(u)\Ai(v)}{u-v},
             \label{airykernellimit2}
        \end{align}
    where $c_2$ is the constant  appearing in
    {\rm(\ref{rhoz2})}.

    Similar results hold near $-z_1$ and $-z_2$.
\end{theorem}

\subsection{The Riemann-Hilbert problem}

Our proofs are based on the fact that the kernel $K_n$ can be
written in terms of
the solution of the following Riemann-Hilbert (RH) problem, see
\cite{DK2}. We look for a
matrix valued function $Y:\mathbb{C}\setminus \mathbb{R}
\to \mathbb{C}^{4\times 4}$ such that

\subsubsection*{RH problem for $Y$}
\begin{enumerate}
    \item[(1)] $Y$ is analytic on $\mathbb{C}\setminus \mathbb{R}$,
    \item[(2)] for $x \in \mathbb{R}$, it holds that
    \begin{align}\label{sprongY2}
        Y_+(x)=Y_-(x)
        \begin{pmatrix}
            1 & 0 & w_{1,1}(x)w_{2,1}(x) & w_{1,1}(x)w_{2,2}(x) \\
            0 & 1 & w_{1,2}(x)w_{2,1}(x) & w_{1,2}(x)w_{2,2}(x) \\
            0 & 0 & 1 & 0 \\
            0 & 0 & 0 & 1
        \end{pmatrix},
    \end{align}
    where $Y_+(x)$ ($Y_-(x)$) denotes the limiting value of $Y(z)$
    as $z$ approaches $x$ from the upper (lower) half plane,
    \item[(3)] as $z \to \infty$, we have that
    \begin{align}\label{Y2inf}
        Y(z) = (I+O(1/z))
        \begin{pmatrix}
            z^{n/2}& 0 & 0 & 0 \\
            0 & z^{n/2} & 0 & 0 \\
            0 & 0 & z^{-n/2} & 0 \\
            0 & 0 & 0 & z^{-n/2}
        \end{pmatrix}.
    \end{align}
\end{enumerate}
In (\ref{sprongY2}) we used the functions
\begin{align}\label{defwC4}
    \begin{array}{ll}
        w_{1,1}(x)=e^{-\frac{n}{2t}(x^2-2ax)}, &
        w_{1,2}(x)=e^{-\frac{n}{2t}(x^2+2ax)},\\
        w_{2,1}(x)=e^{-\frac{n}{2(1-t)}(x^2-2bx)}, &
        w_{2,2}(x)=e^{-\frac{n}{2(1-t)}(x^2+2bx)}.
    \end{array}
\end{align}
This RH problem can be seen as a generalization of the RH problem
for multiple orthogonal polynomials \cite{VAGK}, which in
turn is a generalization of the RH problem for
orthogonal polynomials \cite{FIK}.

The kernel $K_n$ then takes the following form, see \cite{DK2}:
\begin{align} \label{defKinY}
    K_n(x,y)
     & = \frac{1}{2 \pi i(x-y)}
      \begin{pmatrix}
    0  & 0 & w_{2,1}(y) & w_{2,2}(y)
     \end{pmatrix}
     Y_+^{-1}(y) Y_+(x)
     \begin{pmatrix} w_{1,1}(x) \\  w_{1,2}(x)\\ 0\\
     0 \end{pmatrix}.
\end{align}
As explained in \cite{DK2} the solution of the above RH problem is unique
and is built out of polynomials that satisfy certain orthogonality conditions
that are a combination of the conditions satisfied by
multiple Hermite polynomials of type I and II. Therefore they were
called multiple Hermite polynomials of mixed type.
The multiple orthogonal polynomials of mixed type and
their RH problem are related to multi-component KP as shown in \cite{AvMV}.

The formula (\ref{defKinY}) may be seen as a Christoffel-Darboux
formula since it equates the sum (\ref{Knseries}) with the
right-hand side of (\ref{defKinY}) which is of the form
\[ \frac{f_1(x) g_1(y) + f_2(x) g_2(y) + f_3(x) g_3(y) + f_4(x) g_4(y)}{x-y} \]
with only four terms $f_j(x) g_j(y)$ in the numerator. A similar
Christoffel-Darboux formula was derived in \cite{BK2,DK1}
for the case of multiple orthogonal polynomials.

We use the Deift/Zhou steepest descent
method for Riemann-Hilbert problems and apply it to the Riemann-Hilbert
problem stated above. This yields strong and uniform asymptotics for $Y$
as $n \to \infty$ and then also for $K_n$ due to the
formula (\ref{defKinY}).

\subsection{Riemann surface}
As in \cite{ABK,BK2,BK3,KVAW,LW} the asymptotic analysis of the Riemann-Hilbert
problem is based on a suitable Riemann surface that in our case is given
by the equation (\ref{algemeenriemannoppervlak}). The determination of the
equation for this surface is more complicated in this case, since we could not
find an explicit differential equation that is satisfied by the
multiple Hermite polynomials of mixed type. We found the
equation (\ref{algemeenriemannoppervlak})
only after numerical experimentation with Maple.

The equation (\ref{algemeenriemannoppervlak}) has $6$ branch points.
For $ab<1/2$ there are two
critical times $0 < t_{c,1} < t_{c,2} < 1$ such that for $t_{c,1}
< t < t_{c,2}$ four branch points are purely imaginary and two
branch points are real, while for $0 < t < t_{c,1}$ and $t_{c,2} <
t < 1$, four branch points are real and two branch points are
purely imaginary. We prove this in Section \ref{appendix: Branch
points of the Riemann surface}.

The rest of the paper is organized as follows.
In Section \ref{section3} we give more
details about the Riemann surface associated to (\ref{algemeenriemannoppervlak}) that will
be used in the asymptotic analysis of the Riemann-Hilbert problem
for $Y$, defined in (\ref{sprongY2}) and (\ref{Y2inf}). In Section
\ref{section4} we use the
Deift/Zhou steepest descent method to analyze this Riemann-Hilbert
problem. It leads to the proofs of Theorems \ref{theorem21}--\ref{theorem23}
in Section \ref{section5}.

\section{The Riemann surface} \label{section3}

In this section, we study the Riemann surface associated with
the equation (\ref{algemeenriemannoppervlak}) which will be used
in Section \ref{section4} for
the asymptotic analysis of the Riemann-Hilbert problem.
When we talk about the Riemann surface we will always
mean the compact surface that arises after resolution of the singularities
of (\ref{algemeenriemannoppervlak}).

\subsection{Nodal singularities and genus}
First of all we note that (\ref{algemeenriemannoppervlak})
is singular, since clearly both partial derivatives
vanish at $\xi = z = 0$ and so the origin is a singular point.
There are also singular points at infinity and to study
those it is convenient to introduce homogeneous coordinates
$[z:\xi:w]$ and the homogeneous equation
\begin{align} \label{projeq1}
\begin{split}
    & \xi^4-\frac{2}{t(1-t)} z\xi^3+ \frac{1}{t^2(1-t)^2}z^2 \xi^2 +
    \left(-\frac{a^2}{t^2}-\frac{b^2}{(1-t)^2}  + \frac{1}{t(1-t)}\right)\xi^2 w^2
    \\ &+ \left(\frac{2b^2}{t(1-t)^3}-\frac{1}{t^2(1-t)^2}\right)z\xi w^2
    -\frac{b^2}{t^2(1-t)^4} z^2 w^2=0
    \end{split}
\end{align}
in projective space $\mathbb P^2$.

\begin{proposition} \label{prop31}
The equation {\rm (\ref{projeq1})} has three nodal singularities at
$[0:0:w]$, $[z:0:0]$, and $[z: \frac{z}{t(1-t)}:0]$.
Each of these singularities corresponds to two points
on the Riemann surface, and the Riemann surface has
genus $0$.
\end{proposition}
\begin{proof}
To verify that $[z:0:0]$ is a node, we set $z=1$ in (\ref{projeq1})
and we note that the resulting equation in $\xi$ and $w$ has
vanishing partial derivatives at $\xi=w=0$. The quadratic part (in $\xi$ and $w$)
is \[ \frac{1}{t^2(1-t)^2} \xi^2 - \frac{b^2}{t^2(1-t)^4} w^2 \]
which gives rise to two distinct tangents
\begin{equation} \label{tangent}
    \xi = \pm \frac{b}{1-t} w
\end{equation}
and therefore $[z:0:0]$ is a node. Similarly, $[0:0:w]$ is a node.

The third singularity at $[z: \frac{z}{t(1-t)}: 0]$ can be readily seen
if we rewrite (\ref{projeq1}) as
\begin{equation} \label{projeq2}
    \begin{split}
    &
    \xi^2\left(\xi- \frac{z}{t(1-t)}\right)^2 - \frac{a^2}{t^2} \xi^2 w^2
    + \frac{1}{t(1-t)} \xi w^2 \left(\xi-\frac{z}{t(1-t)} \right) \\
    & \qquad -
        \frac{b^2}{(1-t)^2} w^2 \left(\xi - \frac{z}{t(1-t)}\right)^2 = 0
        \end{split}
\end{equation}
Then replacing $z$ by $u = \xi - \frac{z}{t(1-t)}$, setting $\xi=1$
and taking the quadratic part (in $w$ and $u$) we get
$u^2 - \frac{a^2}{t^2} w^2$ which again gives two distinct tangents
$u = \pm \frac{a}{t} w$. Therefore $[z: \frac{z}{t(1-t)}: 0]$
is a node as well.

After resolution  each of the three nodes leads to
two points on the Riemann surface.

The genus formula (Pl\"ucker's formula
\cite[Proposition 2.15]{Mir}) says that
\[ g = (d-1)(d-2)/2 - k \]
for the genus $g$ of a surface of degree $d$ with $k$ nodes
and no other singularities. In the present case we have
$d=4$ and $k \geq 3$. Since $(d-1)(d-2)/2 -k$
is always an upper bound for the genus, we conclude that
the Riemann surface has genus zero. We also conclude that
there are no other singularities besides the three nodal
singularities we already found.
\end{proof}

\subsection{Rational parametrization}
We find a rational parametrization of the Riemann surface by
intersecting the conic
\begin{equation} \label{conic}
    \xi^2 - \frac{1}{t(1-t)} z \xi + p \xi w + q z w = 0
\end{equation}
with the equation (\ref{projeq1}).
By B\'ezout's theorem there are $8$ intersection points in $\mathbb P^2$
if we count intersections according to their multiplicities.
It is easy to see that the conic intersects equation (\ref{projeq1})
at the three nodes $[0:0:w]$, $[z:0:0]$, and $[z: \frac{z}{t(1-t)}: 0]$
(see Proposition \ref{prop31})
for any choice of parameters $p$ and $q$. This accounts for
at least $6$ intersection points, since each of the nodes counts at least twice.
If we choose
\[ q= \frac{b}{t(1-t)^2} \]
then the tangent of (\ref{conic}) at $[z:0:0]$
coincides with one of the tangents (\ref{tangent}). Then we have
higher order intersection at $[z:0:0]$, so that we have to count
this node three times. Then we already have $7$ intersection points.
The remaining intersection point is a point on the surface that is in
one-to-one correspondence with the parameter $p$ and this gives us
the desired parametrization.

Taking $v = t\left(p+ \frac{b}{1-t}\right)$ as a new parameter,
we find after simple calculation that
\begin{align} \label{xiinv}
    \xi & = \frac{b v^2 -v + a^2 b}{(1-t)(a^2 - v^2)} \\
    \nonumber
    z & = t(1-t) \left[ \xi + \frac{v\xi}{t(\xi- \frac{b}{1-t})} \right] \\
      \label{zinv}
      & = \frac{(bv^2 - v + a^2 b)((1-t)v^2 - 2tbv + t -
      (1-t)a^2)}{(2bv-1)(v^2-a^2)}.
\end{align}
The equations (\ref{xiinv}) and (\ref{zinv}) parametrize the surface.
They give a bijection with the Riemann sphere in the $v$-variable
and we conclude once again that the surface has genus zero.

\subsection{The branch points} \label{appendix:
Branch points of the Riemann surface}

Now we start to view (\ref{algemeenriemannoppervlak}) as a
four-sheeted branched covering of the Riemann sphere. That is, for each $z$ the
equation (\ref{algemeenriemannoppervlak}) has four solutions
for $\xi$ where as always we count according to multiplicity. So
each $z$ gives rise to four points or less on the Riemann surface.
The branch points correspond to values of $z$ for which
there are at most three points on the Riemann surface.
The critical $t$-values
\begin{equation} \label{tc1andtc2}
    t_{c,1} = \frac{1+2a^2 - \sqrt{1-4a^2b^2}}{2(1+a^2+b^2)},
    \qquad
    t_{c,2} = \frac{1+2a^2 + \sqrt{1-4a^2b^2}}{2(1+a^2+b^2)}
\end{equation}
are such that
\begin{equation} \label{criticaleq}
     a^2(1-t)^2+ b^2t^2 = t(1-t).
\end{equation}
In that case, when we put $z=0$ in (\ref{algemeenriemannoppervlak}) we obtain
$\xi^4=0$, which means that the two points on the Riemann
surface that correspond to $z=\xi = 0$ are both branch points.
This is the critical case and we will not consider it any further.
It is easy to see that for $ab < 1/2$, we have
$0 < t_{c,1} < t_{c,2} < 1$.

The branch points can be found by calculating the zeros of the
discriminant of (\ref{algemeenriemannoppervlak}). With the aid
of Maple, we get the
following equation for the discriminant:
\begin{equation}\label{discrim}
    (a_3z^6+a_2z^4+a_1z^2+a_0)z^2 = 0,
\end{equation}
where
\begin{align}
    a_3  = & 16a^2b^2, \\
         \nonumber
        a_2  =& (-48a^2b^2(a-b)^2 + a^4 + b^4) t^2 \\
        & + (96a^2(a^2-b^2) - 8 a^2) t - 48a^4b^2-4a^2, \\
        a_1  =& \nonumber
        \left(48a^2b^2(a^4 + b^4 -7a^2b^2-a^2-b^2) \right. \\
        & \qquad \nonumber
            \left. +104a^2b^2 -8a^4-8b^4+20a^2+20b^2+1\right)t^4 \\
        & +\nonumber
            \left(48a^2b^2(4a^4+14a^2b^2+3a^2+b^2) \right. \\
        & \qquad \nonumber
            \left. +32a^4-208a^2b^2-60a^2-20b^2-2\right)t^3 \\
        & + \nonumber (48a^2b^2(6a^4-7a^2b^2 -3a^2) + 104a^2b^2 -48 a^4+60a^2+1)t^2 \\
        & +(48a^4b^2(-4a^2+1)+32a^4-20a^2) t +
        48a^6b^2-8a^4,
        \label{eerstegraadscoefficient} \\
    a_0 =&4(1-4a^2b^2)(a^2(1-t)^2+b^2t^2 - t(1-t))^3.\label{constanteterm}
\end{align}
Note that $z=0$ is always a zero of the discriminant
(\ref{discrim}), but since $z=\xi=0$ is a node, this is not
a branch point in case  $t \neq t_{c,1}$ and $t \neq t_{c,2}$.
So (\ref{discrim}) leads to the sixth degree equation
\begin{align} \label{discrimreduced}
    a_3z^6+a_2z^4+a_1z^2+a_0 = 0,
\end{align}
which has six roots.

We show the following:

\begin{lemma} \label{lemroots}
For every $t \in (0,1)$, the equation {\rm (\ref{discrimreduced})}
has only real or purely imaginary roots.
\begin{enumerate}
\item[\rm (a)] For $t \in (0,t_{c,1}) \cup (t_{c,2},0)$, there are four real roots
$\pm z_1$ and $\pm z_2$ with $z_1 > z_2 > 0$ and two imaginary
roots $\pm i z_3$ with $z_3 > 0$.
\item[\rm (b)] For $t \in (t_{c,1}, t_{c,2})$, there are two
real roots $\pm z_1$ with $z_1 > 0$ and four purely imaginary
roots $\pm i z_2$ and $\pm iz_3$ with $z_3 \geq z_2 > 0$.
\end{enumerate}
\end{lemma}

\begin{proof}
 Because (\ref{discrimreduced}) only has even powers of $z$,
 it is enough to prove that the polynomial
\begin{equation}\label{veelterm}
    p_1(x) = a_3x^3+a_2x^2+a_1x+a_0
\end{equation}
has three real roots with one negative root and two positive
roots in case (a) and one positive root and two negative roots
in case (b).

We first examine when (\ref{veelterm}) has a multiple root.
We find those values by looking at the zeros of the discriminant
of (\ref{veelterm}), which is a polynomial of degree $12$ in $t$,
given by
    \begin{equation}
        \begin{split}
        &
            16t(1-t)((a+b)t-a)^2((a-b)t-a)^2\\ & \times\left((432a^4b^4+8a^2-72a^2b^2-1+8b^2)t^2\right.
            \\
            & \qquad \left.+((-432a^4b^4+ 72a^2b^2+1)t(1-t)+8a^2(1-t)^2+8b^2t^2\right)^3,
        \end{split}
    \end{equation}
which was calculated with the use of Maple.
We find multiple roots for $t_1=0$, $t_2=1$, $t_3=\frac{a}{a+b}$,
$t_4 = \frac{a}{a-b}$, and for the zeros $t_5$ and $t_6$
of the quadratic polynomial
\[ q(t) = (-432a^4b^4+ 72a^2b^2+1)t(1-t)+8a^2(1-t)^2+8b^2t^2. \]
Now it is easy to verify that
\[ 8 a^2 \frac{1-t}{t} + 8 b^2\frac{t}{1-t} \]
takes its minimum for $t \in (0,1)$ in $t= t_3 =\frac{a}{a+b}$
and the minimum is $16ab$.
Thus for $t \in (0,1)$,
\begin{align*}
    q(t) & \geq t(1-t) (-432a^4b^4+ 72a^2b^2+1 + 16ab) \\
    & =  t(1-t) (1-2ab)(1+6ab)^3 > 0,
 \end{align*}
since $ab < 1/2$. So $t_5$ and $t_6$ do not belong to
the interval $(0,1)$ and the only value of $t \in (0,1)$
for which (\ref{veelterm}) has a multiple root is
$t_3 = \frac{a}{a+b}$.
For $t_3$ we find that $p_1$ has the following roots:
\begin{align*}
        x_1 = \frac{4ab(2ab+1)}{(a+b)^2} \qquad \mbox{ and }
        \qquad x_2 = x_3 = -\frac{(2ab-1)^2}{4(a+b)^2}.
    \end{align*}
So for $t=t_3$ we have that $p_1$ has a double negative
root and one positive root.

Next, we investigate the critical $t$-values
$t = t_{c,1}$ and $t=t_{c,2}$, see (\ref{tc1andtc2}),
which are exactly those $t \in (0,1)$ for which
$a_0 = 0$.
So for these values, one root of (\ref{veelterm}) is at $x=0$.
We determine the sign of the other two roots by looking
at $a_1$. We are going to take $t=t_{c,1}$. The proof for $t=t_{c,2}$
is similar. If we substitute $t=t_{c,1}$ into
(\ref{eerstegraadscoefficient}), we arrive at the expression
\begin{align} \nonumber
    a_1 &= -\frac{27(1-4a^2b^2)^2}{2(a^2+b^2+1)^4}  \left(
        a^4 ((4b^2+1)(1+\sqrt{1-4a^2b^2}) -2a^2b^2) \right. \\
        & \qquad \qquad \qquad \label{a1atcriticaltime}
        \left. + b^4 ((4a^2+1)(1-\sqrt{1-4a^2b^2}) -2a^2b^2) + 12a^4 b^4 \right)
\end{align}
Since $0 < ab < 1/2$, it is easy to check that  $2a^2b^2 < 1-\sqrt{1-4a^2b^2}$.
Then it readily follows that $a_1 < 0$ because of (\ref{a1atcriticaltime}).
Thus for $t=t_{c,1}$ the polynomial $p_1$ has one negative root,
one positive root, and one root at $0$. The same holds for $t=t_{c,2}$.

Now we can finish the proof of the lemma by a continuity argument.
For $t=t_{c,1}$ and for $t=t_{c,2}$ we have three distinct real roots.
Since there are only multiple roots for $t=a/(a+b)$ we find three
three distinct real roots for every $t$ in the two
intervals $(0,a/(a+b))$ and $(a/(a+b),1)$.
Note that $t_{c,1} < a/(a+b) < t_{c,2}$.

Since we have a double negative root and one positive root for $t=a/(a+b)$,
and $0$ is only a root for $t_{c,1}$ and $t_{c,2}$, it follows
that $p_1$ has two negative roots and one positive root for
every $t \in (t_{c,1}, t_{c,2})$.

For $t=0$, we find $p_1(x) = 4a^2(x-a^2)^2(4b^2 x - 4a^2b^2+1)$
with roots $x_1 = x_2 = a^2 > 0$ and $x_3 = \frac{4a^2b^2-1}{4b^2} < 0$.
Then by continuity we find that there are two positive roots and
one negative root  for every $t \in (0,t_{c,1})$. We find
the same thing for $t \in (t_{c,2},1)$ by looking at $t=1$.

This completes the proof of the lemma.
\end{proof}

The six roots of (\ref{discrimreduced}) are all simple
branch points of the Riemann surface. This follows from the
Riemann-Hurwitz formula \cite[Theorem 4.16]{Mir}, and the fact that the genus is zero.

\subsection{Sheet structure}

There are four inverse functions of (\ref{algemeenriemannoppervlak}), which
behave near infinity as
\begin{align}
    \xi_1(z) & = \frac{z}{t(1-t)}-\frac{a}{t}-\frac{1}{2z} + O\left(\frac{1}{z^2}\right),\label{infxi21} \\
    \xi_2(z) & = \frac{z}{t(1-t)}+\frac{a}{t}-\frac{1}{2z} + O\left(\frac{1}{z^2}\right),\label{infxi22} \\
    \xi_3(z) & = \frac{b}{1-t}+\frac{1}{2z} + O\left(\frac{1}{z^2}\right), \label{infxi23} \\
    \xi_4(z) & = -\frac{b}{1-t} + \frac{1}{2z} + O\left(\frac{1}{z^2}\right). \label{infxi24}
\end{align}

\begin{figure}
    \begin{center}
        \includegraphics[height=7cm,width=10cm]{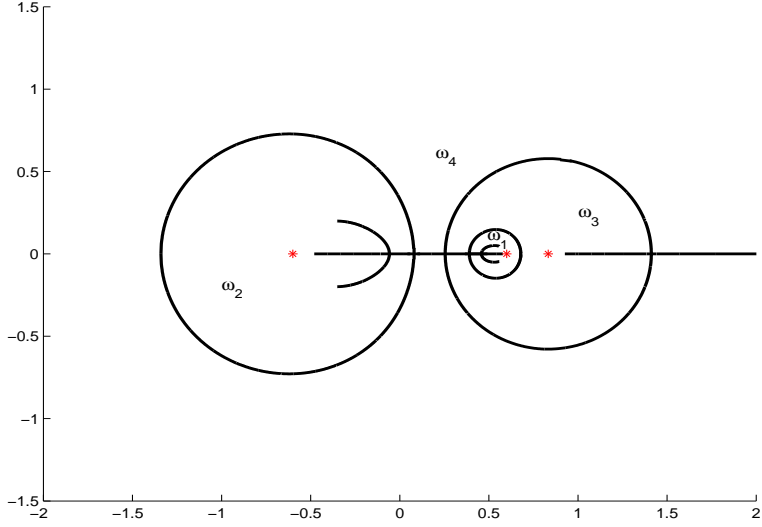}
        \caption{\label{vplanet=0.25}
       The image of the cross $C = [-z_1,z_1] \cup [-iz_3,iz_3]$
       in the $v$-plane for   $0 < t < t_{c,1}$. The values
       $v=-a$, $v=a$, and $v= \frac{1}{2b}$ are indicated with $\ast$.
       These $v$-values correspond to $\infty$ on the second, first
       and third sheets, respectively. The precise values used
       for the figure are $a=b=0.6$ and $t=0.25$.}
    \end{center}
\end{figure}

The sheet structure of the Riemann surface is determined by the
way we choose the analytical continuations of the $\xi$-functions.
Since all branch points are on the real or imaginary axis,
we find that all $\xi_j$'s have an analytic extension to
$\mathbb C \setminus C$ where $C$ is the cross
\[ C = [-z_1,z_1] \cup [-iz_3, iz_3]. \]
The analytic continuation is also denoted by $\xi_j$.
We consider $\xi_j$ on the $j$th sheet.

\begin{figure}
    \begin{center}
        \includegraphics[height=7cm,width=10cm]{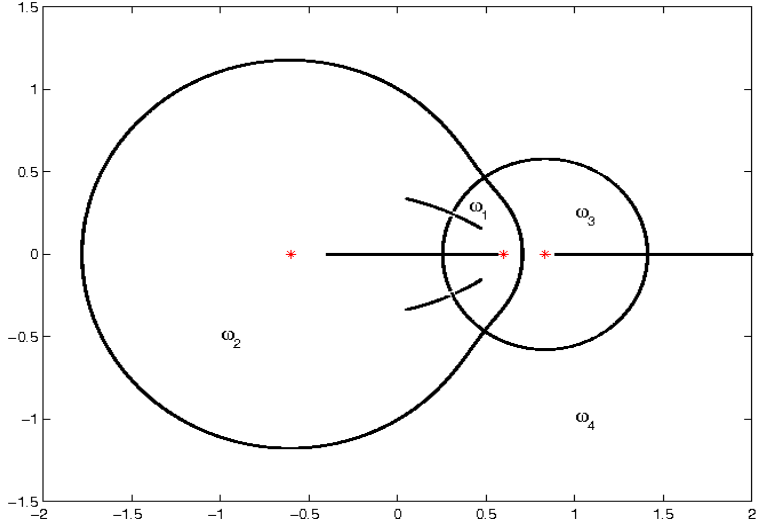}
        \caption{\label{vplanet=0.45}
        The image of the cross $C = [-z_1,z_1] \cup [-iz_3,iz_3]$ in the $v$-plane for
        $t_{c,1}< t < t_{c,2}$. The values
       $v=-a$, $v=a$, and $v= \frac{1}{2b}$ are indicated with $\ast$.
       These $v$-values correspond to $\infty$ on the second, first
       and third sheets, respectively. The precise values used
       for the figure are $a=b=0.6$ and $t=0.45$.}
    \end{center}
\end{figure}

Now we use the rational parametrization (\ref{xiinv})--(\ref{zinv})
of the surface to compute the corresponding curves in the
$v$-plane. We denote the $v$-value corresponding
to $z \in \mathbb C \setminus C$ on the $j$th sheet
by $v_j(z)$ and we let $\omega_j$ be the image set
of $v_j(\mathbb C \setminus C)$.
The next two figures show the regions
$\omega_j$, $j=1,2,3,4$, for the cases $t < t_{c,1}$
and $t_{c,1} < t < t_{c,2}$, respectively.
Denoting the point at infinity on the $j$th sheet
by $\infty_j$, we can determine the corresponding
$v$-values from the formulas (\ref{xiinv}), (\ref{zinv}),
and (\ref{infxi21})--(\ref{infxi24}). Indeed we have
\begin{align} \label{vatinfty}
    \left\{  \begin{array}{ll}
    z = \infty_1 \ \longleftrightarrow \ v = a, \\
    z = \infty_2 \ \longleftrightarrow \ v = -a, \\
    z = \infty_3 \ \longleftrightarrow \ v = \frac{1}{2b},  \\
    z = \infty_4 \ \longleftrightarrow \ v = \infty.
    \end{array} \right.
\end{align}
These four $v$-values are independent of $t$, and they
are indicated with $\ast$ in Figures \ref{vplanet=0.25}
and \ref{vplanet=0.45} (except for $v=\infty$, which is of
course not visible). Note that $-a < a < \frac{1}{2b}$.

The figures show what the analytic continuations
are of the functions $\xi_j$ and how the sheets of the
Riemann surface are connected.
First of all, we note that any part of the boundary
of $\omega_j$ that is not a part of the boundary of some
$\omega_k$ with $k\neq j$, corresponds in the $z$-plane
to a part of the cross $C$ where $\xi_j$ has equal boundary
values from both sides. This part of the cross is then
removed from the cut on the $j$th sheet.

\begin{figure}[t]
    \begin{center}
        \includegraphics[height=7cm]{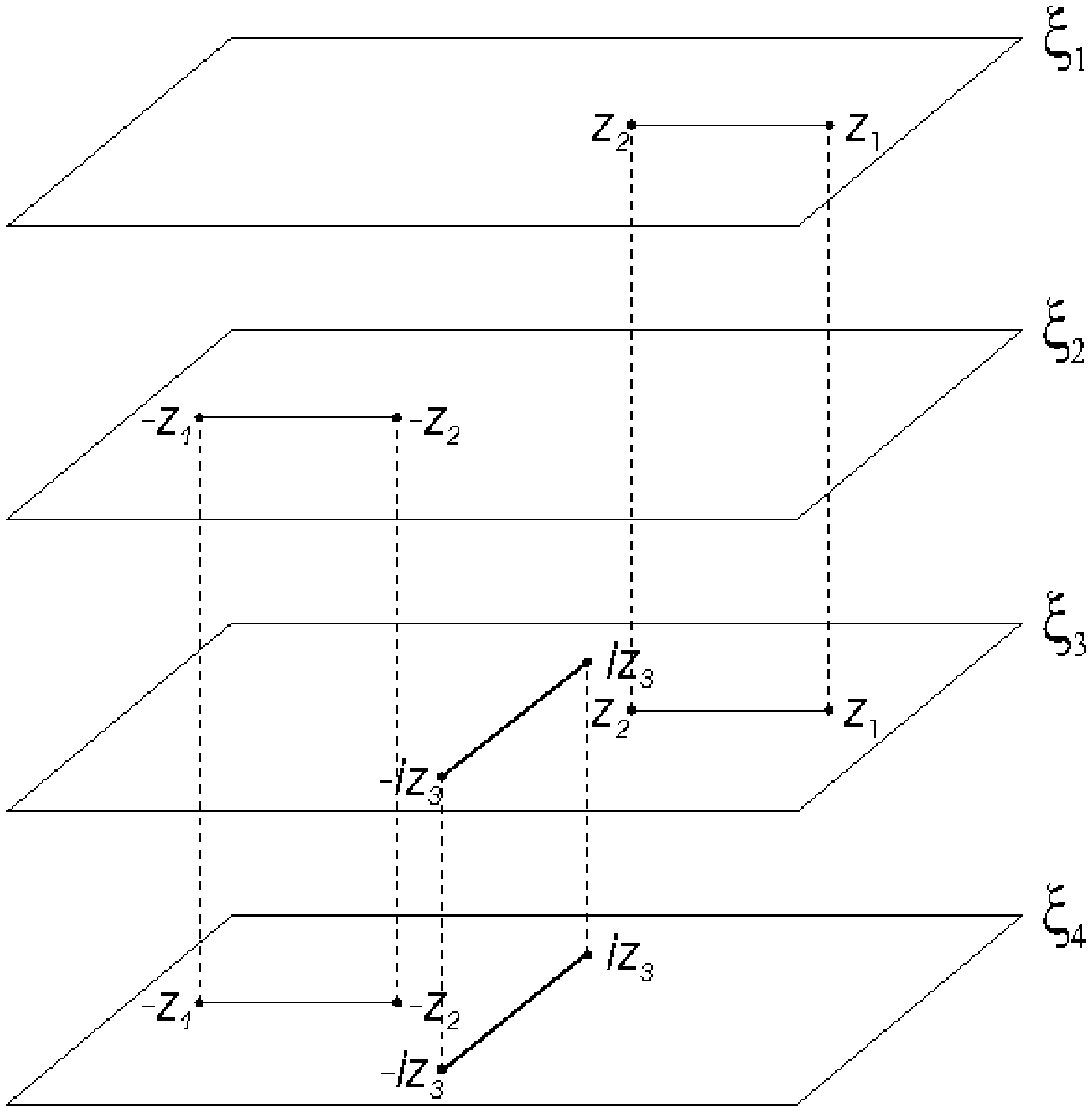}
        \caption{\label{5figuur1}
        The sheet structure of the Riemann surface for $ab<1/2$ and $0 < t < t_{c,1}$.}
    \end{center}
\end{figure}

In this way, we see that for $t < t_{c,1}$,
$\xi_1$ is defined with a cut $[z_2,z_1]$ only,
$\xi_2$ with a cut $[-z_1,-z_2]$,
$\xi_3$ with two cuts $[z_2,z_1]$ and $[-iz_3, iz_3]$
and $\xi_4$ is defined with two cuts
$[-z_2, -z_1]$ and $[-iz_3,iz_3]$.
The sheets are connected as shown in Figure~\ref{5figuur1}.

\begin{figure}
    \begin{center}
        \includegraphics[height=8cm]{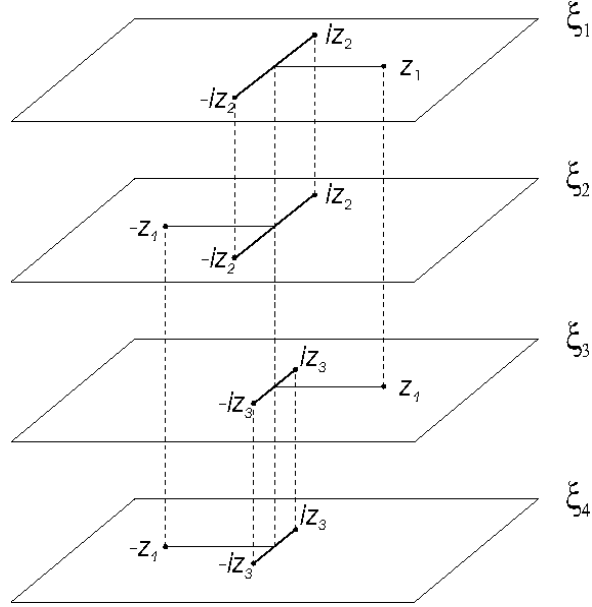}
        \caption{\label{figuur9f}
        The sheet structure of the Riemann surface for $ab<1/2$ and
        $t \in (t_{c,1},t_{c,2})$.}
    \end{center}
\end{figure}

For $t_{c,1} < t < t_{c,2}$ we have the two pairs of purely
imaginary branch points $\pm iz_2$ and $\pm iz_3$,
where we choose $z_3 \geq z_2$.
In the situation of Figure~\ref{vplanet=0.45}
we have that $\xi_1(i z_2) = \xi_2(iz_2)$ and $\xi_3(iz_3) = \xi_4(iz_3)$.
This is always the case if $t_{c,1} < t < a/(a+b)$.
Then we have that $\xi_1$ is analytic with a
cut on $[0,z_1] \cup [-iz_2,iz_2]$, $\xi_2$
has a cut on $[-z_1,0] \cup [-iz_2,iz_2]$,
$\xi_3$ has a cut on $[0,z_1] \cup [-iz_3,iz_3]$,
and $\xi_4$ has a cut on $[-z_1,0] \cup [-iz_3,iz_3]$.
The sheet structure is then as shown in Figure~\ref{figuur9f}.

For $t=a/(a+b)$ we have that $z_2 = z_3$ and
for $a/(a+b) < t < t_{c,2}$ we find that the
role of $z_2$ and $z_3$ are reversed, so that
then the first and second sheets are
connected along $[-iz_3,iz_3]$ and the third
and fourth along $[-iz_2,iz_2]$.

For $t_{c,2} < t < 1$, we have a similar sheet structure as in
Figure~\ref{5figuur1}, except that the cut on the vertical segment
$[-iz_3,iz_3]$ now connects the first and the second sheets.

\subsection{Properties of $\xi_j$}

The sheet structure induces jump relations between the
$\xi$-functions along the cuts that we will use later on.
At the branch point $z_1$ we have for a real
constant $c_1 > 0$ that
\begin{align}
    \xi_1(z) &= \xi_1(z_1)+c_1(z-z_1)^{1/2}+O(z-z_1), \label{rhorondz1}
    \\
    \xi_3(z) &= \xi_1(z_1)-c_1(z-z_1)^{1/2}+O(z-z_1).
\end{align}
In case $t \in (0,t_{c,1}) \cup (t_{c,2},1)$ we also have
that at the branch point $z_2$ there is a real constant
$c_2>0$ such that
\begin{align}
    \xi_1(z) &= \xi_1(z_2)-c_2(z_2-z)^{1/2}+O(z_2-z),\label{rhorond-z1}\\
    \xi_3(z) & = \xi_1(z_2) + c_2(z_2-z)^{1/2}+O(z_2-z).
\end{align}
Here the square root is taken with a branch cut along the negative
real axis. Similar behavior also holds for the functions $\xi_2$
and $\xi_4$ near the branch points $-z_1$ and $-z_2$.

From the sheet structure we also see how the
$\xi$-functions are continued along the cuts.
We will not list those relations explicitly,
but they will be (tacitly) used in the future.

What is also important is that for any closed contour $\gamma$
that does not intersect the cut on the $j$th sheet, we have
\begin{equation} \label{integerloops}
    \frac{1}{2\pi i} \oint_{\gamma} \xi_j(s) ds \in
    \frac{1}{2} \mathbb Z.
\end{equation}
This follows from the fact that the residues
at infinity for each of the $\xi_j$ functions is $\pm \frac{1}{2}$,
see (\ref{infxi21})--(\ref{infxi24}).

\subsection{The $\lambda_j$ functions}
The $\lambda$-functions are primitive functions of
the $\xi$-functions. We define them as
\begin{align}\label{deflambda2}
    \begin{split}
    \lambda_1(z) & =  \int_{z_1}^z \xi_1(s)ds,\\
    \lambda_2(z) &=  \int_{-z_{1+}}^z \xi_2(s)ds + c,\\
    \lambda_3(z) &=  \int_{z_1}^z \xi_3(s)ds, \\
    \lambda_4(z) &=  \int_{-z_{1-}}^z \xi_4(s)ds + c,
    \end{split}
\end{align}
where $c$ is chosen so that $\lambda_3(iz_3) = \lambda_4(iz_3)$
(if $0 < t < a/(a+b)$), or $\lambda_1(iz_3) = \lambda_2(iz_3)$
(if $a/(a+b) < t < 1$).
The path of integration in each integral in (\ref{deflambda2})
is in $\mathbb C \setminus (C \cup (-\infty,-z_1])$.

From these definitions, and from the fact (\ref{integerloops})
that the periods of the $\xi$-functions are half integers,
it then follows that $e^{n\lambda}$ is analytic on the
Riemann surface (recall $n$ is even). That is, if $C_{jk}$
is the cut connecting the $j$th and $k$th sheets, then
\begin{equation} \label{continuationlambda}
    e^{n\lambda_{j\pm}} = e^{n\lambda_{k\mp}} \qquad \mbox{on } C_{jk}.
\end{equation}

The behavior as $z\to \infty$ follows from (\ref{deflambda2}) and
(\ref{infxi21})--(\ref{infxi24}) and is given by
\begin{align}\label{lambda1inf}
    \lambda_1(z) & =  \frac{z^2}{2t(1-t)} -\frac{az}{t} -\frac{1}{2}\ln z +l_1 +O\left(\frac{1}{z}\right), \\
    \label{lambda2inf}
    \lambda_2(z) & =  \frac{z^2}{2t(1-t)} +\frac{az}{t} -\frac{1}{2}\ln z +l_2 +O\left(\frac{1}{z}\right), \\
    \label{lambda3inf}
    \lambda_3(z) & =  \frac{bz}{1-t} + \frac{1}{2}\ln z +l_3 +O\left(\frac{1}{z}\right), \\
    \lambda_4(z) & =  -\frac{bz}{1-t} + \frac{1}{2}\ln z +l_4
    +O\left(\frac{1}{z}\right), \label{lambda4inf}
\end{align}
where $l_1,l_2,l_3,l_4$ are certain integration constants.

\section{Steepest descent analysis for case $0 < t < t_{c,1}$}
\label{section4}

We will do the steepest descent analysis in some detail for the
case $0 < t < t_{c,1}$. The case $t_{c,2} < t < 1$ follows
from this by symmetry of the problem. For the case $t_{c,1} < t < t_{c,2}$,
we refer to Section \ref{section6}.

The steepest descent analysis is based on the one given in \cite{ABK}.
A major role is played by the Riemann surface, which for
for $0 < t < t_{c,1}$, has the sheet structure as shown
in Figure~\ref{5figuur1}.

\subsection{First transformation: $Y \mapsto U$}
In the first transformation we normalize the RH problem
at infinity. We define
\begin{align}\label{defU5}
    \begin{split}
    U(z)& = L^n Y(z)\diag\left(e^{n(\lambda_1(z)-\frac{z^2}{2t(1-t)}+\frac{az}{t})},e^{n(\lambda_2(z)-\frac{z^2}{2t(1-t)}-\frac{az}{t})}
    ,\right.\\& \qquad \qquad \left.e^{n(\lambda_3(z)-\frac{bz}{1-t})}, e^{n(\lambda_4(z)+\frac{bz}{1-t})}\right),
    \end{split}
\end{align}
where $L$ is the constant diagonal matrix
\begin{align} \label{defL}
    L & = \diag\left(e^{-l_1}, e^{-l_2},
        e^{-l_3}, e^{-l_4}\right).
\end{align}
Here the constants $l_j$, $j=1,2,3,4$, are the constants
that appear in (\ref{lambda1inf})--(\ref{lambda4inf}).
Now $U$ is  defined and analytic in
$\mathbb{C} \setminus (\mathbb{R} \cup [-iz_3,iz_3])$. Then
(\ref{Y2inf}), (\ref{lambda1inf})--(\ref{lambda4inf}),
and (\ref{defU5}) imply that $U$ is normalized at infinity:
\begin{align}
    U(z) = I+ O(1/z) \quad \mbox{ as $n \to \infty$}.
\end{align}
 Using (\ref{defU5}) and (\ref{sprongY2}) we find the jumps for
 $U$. On the real line, we get
 \begin{align} \label{jumpU0}
    U_+ &=U_-
    \begin{pmatrix}
        e^{n(\lambda_{1+}-\lambda_{1-})} & 0 & e^{n(\lambda_{3+}-\lambda_{1-})} & e^{n(\lambda_{4+}-\lambda_{1-})} \\
        0 & e^{n(\lambda_{2+}-\lambda_{2-})} &
        e^{n(\lambda_{3+}-\lambda_{2-})}&
        e^{n(\lambda_{4+}-\lambda_{2-})}\\
        0 & 0 & e^{n(\lambda_{3+}-\lambda_{3-})} & 0 \\
        0 & 0 & 0 & e^{n(\lambda_{4+}-\lambda_{4-})}
    \end{pmatrix},
 \end{align}
 while on the vertical segment $[-iz_3,iz_3]$, which is oriented upwards (so that
 $U_+(z)$ ($U_-(z))$ is the limiting value of $U(z')$ as $z' \to z$ from the left (right) half-plane)
  we have the jump
 \begin{align}
    U_+ & = U_-
        \diag\left(1, 1, e^{n(\lambda_{3+}-\lambda_{3-})},
        e^{n(\lambda_{4+}-\lambda_{4-})}\right).
 \end{align}
Because of (\ref{continuationlambda}) we can simplify the jump (\ref{jumpU0}) on the
various parts of the real line. The result is the following RH problem for $U$:

\subsubsection*{RH problem for $U$}
\begin{enumerate}
    \item[(1)] $U$ is analytic on $\mathbb{C} \setminus (\mathbb{R}
     \cup [-iz_3,iz_3])$.
    \item[(2)] $U$ satisfies the following jumps on the real line:
    \begin{align}
        U_+ & = U_-
        \begin{pmatrix}
            1 & 0 & e^{n(\lambda_{3+}-\lambda_{1-})} &
            e^{n(\lambda_{4+}-\lambda_{1-})} \\
            0 & 1 & e^{n(\lambda_{3+}-\lambda_{2-})} &
            e^{n(\lambda_{4+}-\lambda_{2-})} \\
            0 & 0 & 1 & 0 \\
            0 & 0 & 0 & 1
        \end{pmatrix}
        \quad \mbox{on $(-\infty,-z_1)$},\label{jumpU1} \\
        U_+ & = U_-
        \begin{pmatrix}
        1 & 0 & e^{n(\lambda_{3+}-\lambda_{1-})} &
        e^{n(\lambda_{4+}-\lambda_{1-})} \\
        0 & e^{n(\lambda_{2+}-\lambda_{2-})} & e^{n(\lambda_{3+}-\lambda_{2-})} &
        1 \\
        0 & 0 & 1 & 0 \\
        0 & 0 & 0 & e^{n(\lambda_{4+}-\lambda_{4-})}
        \end{pmatrix}
        \quad \mbox{on $(-z_1,-z_2)$}, \label{jumpU2} \\
        U_+ & = U_-
        \begin{pmatrix}
            1 & 0 & e^{n(\lambda_{3+}-\lambda_{1-})} &
            e^{n(\lambda_{4+}-\lambda_{1-})} \\
            0 & 1 & e^{n(\lambda_{3+}-\lambda_{2})} &
            e^{n(\lambda_{4+}-\lambda_{2})} \\
            0 & 0 & 1 & 0 \\
            0 & 0 & 0 & 1
        \end{pmatrix}
        \quad \mbox{on $(-z_2,0) \cup (0,z_2)$}, \label{jumpU3} \\
        U_+ & = U_-
        \begin{pmatrix}
            e^{n(\lambda_{1+}-\lambda_{1-})} & 0 &
            1 &
            e^{n(\lambda_{4}-\lambda_{1-})} \\
            0 & 1 & e^{n(\lambda_{3+}-\lambda_{2})} &
            e^{n(\lambda_{4}-\lambda_{2})}\\
            0 & 0 & e^{n(\lambda_{3+}-\lambda_{3-})} & 0 \\
            0 & 0 & 0 & 1
        \end{pmatrix}
        \quad \mbox{on $(z_2,z_1)$}, \label{jumpU4}
    \end{align}
    \begin{align}
        U_+ & = U_-
        \begin{pmatrix}
            1 & 0 & e^{n(\lambda_{3}-\lambda_{1})} &
            e^{n(\lambda_{4}-\lambda_{1})} \\
            0 & 1 & e^{n(\lambda_{3}-\lambda_{2})} &
            e^{n(\lambda_{4}-\lambda_{2})} \\
            0 & 0 & 1 & 0 \\
            0 & 0 & 0 & 1
        \end{pmatrix}\label{jumpU5}
        \quad \mbox{on $(z_1,+\infty)$}.
    \end{align}
    Finally, on the vertical segment, which is oriented upwards,
    \begin{align}
        U_+ & = U_-\diag\left(1, 1, e^{n(\lambda_{3+}-\lambda_{3-})},
        e^{n(\lambda_{4+}-\lambda_{4-})}\right)
        \quad \mbox{on $(-iz_3, i z_3)$}.
        \label{jumpU6}
    \end{align}
    \item[(3)] $U(z) = I+ O(1/z)$ as $z \to \infty$.
\end{enumerate}

Now it turns out that not all entries in the jump matrices
for $U$ are well-behaved as $n \to \infty$. Ideally, one
would like to have exponentially decaying off-diagonal terms
$e^{n(\lambda_{j+}-\lambda_{k-})}$ in all of the jump matrices
in (\ref{jumpU1})--(\ref{jumpU6}), and oscillating diagonal
entries $e^{n(\lambda_{j+}-\lambda_{j-})}$. However, this is not
the case in the present situation. For example, the entry $e^{n(\lambda_{3+}-\lambda_{3-})}$
in the jump matrix in (\ref{jumpU6}) is exponentially
increasing as $n \to \infty$, since $\Re \lambda_{3+} > \Re\lambda_{3-}$
on the vertical segment $(-iz_3,iz_3)$.

\subsection{Second transformation: $U \mapsto T$}

We use a trick to remove the exponentially increasing
entries. This involves the global opening of a lens, which
was introduced in \cite{ABK}, and also used in \cite{LW}.

We take a closed curve $\Sigma$, consisting of a part in the right
half-plane that connects $-iz_3$ and $iz_3$, that is symmetric
with respect to the real axis, and that satisfies
\begin{equation} \label{orderingoflambdas1}
     \Re \lambda_4(z) < \Re \lambda_3(z)
    \qquad \mbox{ for } z \in \Sigma \cap \{ \Re z > 0\}.
    \end{equation}
The curve $\Sigma$ intersects the positive real line in a point $x^*$
and we assume that $x^*$ is sufficiently large so that
\[ \Re \lambda_4(x) < \Re \lambda_3(x) < \Re \lambda_1(x) < \Re \lambda_2(x)
    \qquad \mbox{ for } x \geq x^*, \ x \in \mathbb R.
    \]
The part of $\Sigma$ in the left half-plane is the mirror image
with respect to the imaginary axis. Then
\begin{equation} \label{orderingoflambdas2}
    \Re \lambda_3(z) <  \Re \lambda_4(z)
    \qquad \mbox{ for } z \in \Sigma \cap \{ \Re z < 0\},
    \end{equation}
and
\[ \Re \lambda_3(x) < \Re \lambda_4(x) < \Re \lambda_2(x) < \Re \lambda_1(x)
    \qquad \mbox{ for } x \leq -x^*, \ x \in \mathbb R.
    \]
Note that the behavior (\ref{orderingoflambdas1})--(\ref{orderingoflambdas2})
is valid near infinity because of (\ref{lambda3inf})--(\ref{lambda4inf}) and
we have that the inequalities remain valid in a domain whose
boundary contains the branch points $\pm i z_3$.

\begin{figure}
    \begin{center}
        \includegraphics[height=8cm]{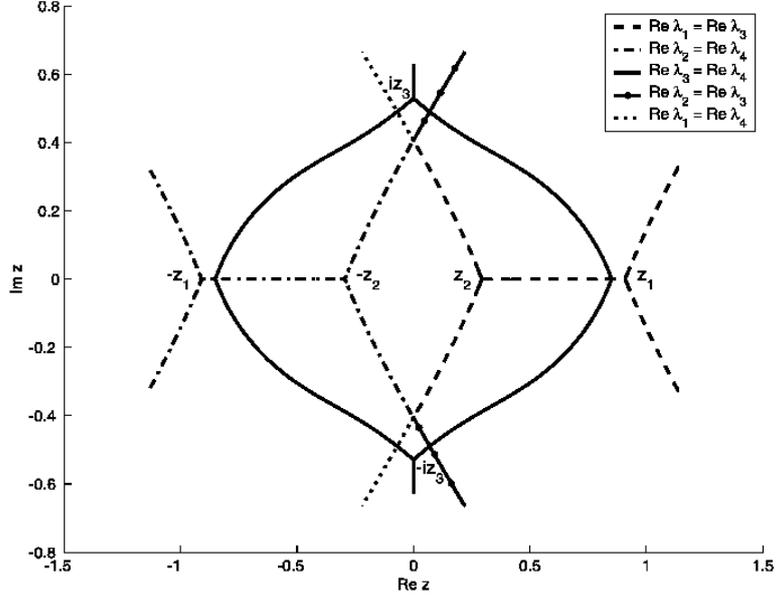}
        \caption{\label{5figuur3a}
        The behavior of the real part of $\lambda_j-\lambda_k$ when $a=b=0.6$
        and $t=0.05$. Only some relevant curves are shown
        where two real parts coincide.}
    \end{center}
\end{figure}

\begin{figure}
    \begin{center}
        \includegraphics[height=7.8cm]{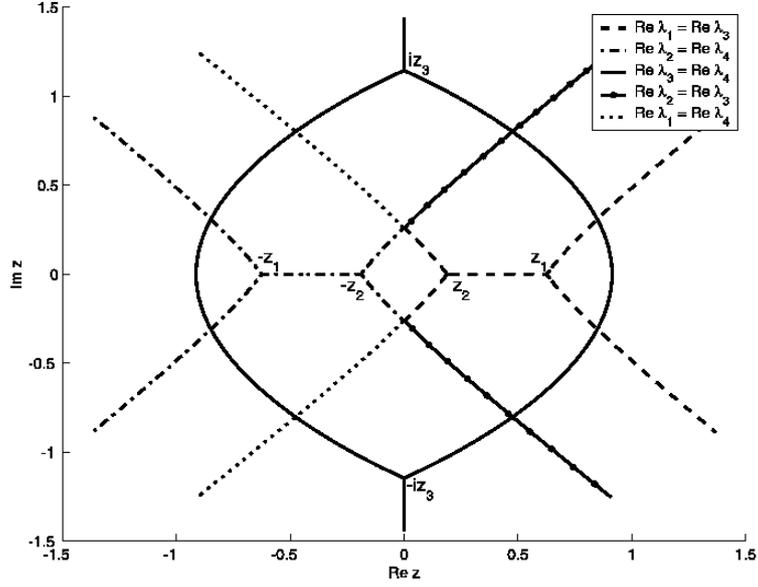}
        \caption{\label{5figuur3b}
        The behavior of the real part of $\lambda_j-\lambda_k$
        when $a=b=0.4$ and $t=0.025$.}
    \end{center}
\end{figure}

This is illustrated in Figures \ref{5figuur3a} and \ref{5figuur3b}
where the solid curves are the curves where
$\Re \lambda_3 = \Re \lambda_4$. We have three such curves that
emanate at equal angles from $i z_3$. One of them extends
to infinity along the positive imaginary axis. The other two continue
to the branch point $-iz_3$ and together form a closed contour. This
closed contour can either contain the intervals $[-z_1,-z_2]$
and $[z_2,z_1]$ in its interior as in Figure~\ref{5figuur3b},
or it intersects these two intervals as in Figure~\ref{5figuur3a}.

\begin{figure}[t]
    \begin{center}
        \includegraphics[height=7.8cm]{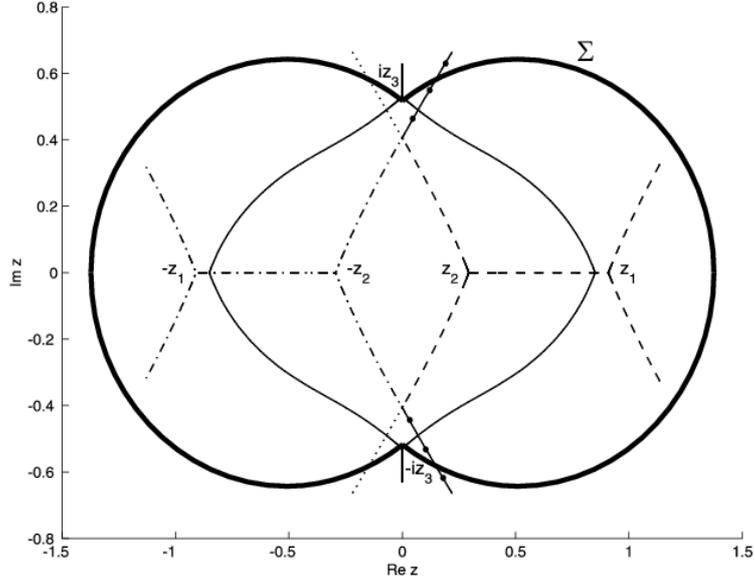}
        \caption{\label{5figuur4} Opening of the global lens $\Sigma$
        that connects $\pm iz_3$. It is such that $\Re \lambda_3 > \Re \lambda_4$
        on $\Sigma$ in the right half-plane, and $\Re \lambda_4 > \Re \lambda_3$
        on $\Sigma$ in the left half-plane. Here we have taken
        $a=b=0.6$ and $t=0.05$.}
    \end{center}
\end{figure}

For later convenience, we also choose $\Sigma$ as the analytic continuation
of the closed contour where $\Re \lambda_3 = \Re \lambda_4$ in
a neighborhood of $\pm iz_3$. For an example how to
choose $\Sigma$, see Figure~\ref{5figuur4}.

Now we are ready to introduce the second transformation,
see also \cite[section 4]{ABK}. We define $T$ as
\begin{align}
    T & = U  \quad \mbox{ outside $\Sigma$}, \\
    T & = U
    \begin{pmatrix}
        1 & 0 & 0 & 0 \\
        0 & 1 & 0 & 0 \\
        0 & 0 & 1 & 0 \\
        0 & 0 & -e^{n(\lambda_3-\lambda_4)} & 1
    \end{pmatrix}
    \quad \mbox{ inside $\Sigma$ in the left half-plane},
\end{align}
\begin{align}
    T & = U
    \begin{pmatrix}
        1 & 0 & 0 & 0 \\
        0 & 1 & 0 & 0 \\
        0 & 0 & 1 & -e^{n(\lambda_4-\lambda_3)} \\
        0 & 0 & 0 & 1
    \end{pmatrix}
    \quad \mbox{ inside $\Sigma$ in the right half-plane}.
    \label{defTbinnenSigmarechts}
\end{align}

Then  $T$ has jumps on the contours shown in Figure~\ref{5figuur5} and a straightforward calculation shows
that $T$ satisfies the following RH problem:

\begin{figure}
    \begin{center}
        \includegraphics[height=5cm]{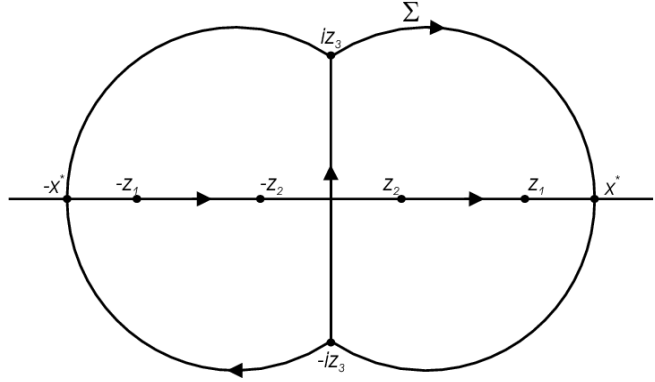}
        \caption{\label{5figuur5}
        Jump contour for the RH problem for $T$. The contour $\Sigma$ has clockwise
        orientation. The $4\times 4$ matrix valued
        function $T$ is analytic outside these contours, and has jumps (\ref{jumpT1})--(\ref{jumpT11})
        along the various parts of the contour. }
    \end{center}
\end{figure}

\subsubsection*{RH problem for $T$}
\begin{enumerate}
    \item[(1)] $T$ is analytic on $\mathbb{C} \setminus(\mathbb{R}
     \cup [-iz_3,iz_3] \cup \Sigma)$.
    \item[(2)] On the real line, $T$ has the following jumps:
    \begin{align}
        T_+ & = T_-
        \begin{pmatrix}
            1 & 0 & e^{n(\lambda_{3+}-\lambda_{1-})} &
            e^{n(\lambda_{4+}-\lambda_{1-})} \\
            0 & 1 & e^{n(\lambda_{3+}-\lambda_{2-})} &
            e^{n(\lambda_{4+}-\lambda_{2-})} \\
            0 & 0 & 1 & 0 \\
            0 & 0 & 0 & 1
        \end{pmatrix}
        \quad \mbox{on $(-\infty, -x^*)$}, \label{jumpT1} \\
        T_+ & = T_-
        \begin{pmatrix}
            1 & 0 & 0 & e^{n(\lambda_{4+}-\lambda_{1-})} \\
            0 & 1 & 0 & e^{n(\lambda_{4+}-\lambda_{2-})} \\
            0 & 0 & 1 & 0 \\
            0 & 0 & 0 & 1
        \end{pmatrix}
        \quad \mbox {on $(-x^*,-z_1)$}, \label{jumpT2}
    \end{align}
    \begin{align}
        T_+ & = T_-
        \begin{pmatrix}
            1 & 0 & 0 & e^{n(\lambda_{4+}-\lambda_{1-})} \\
            0 & e^{n(\lambda_{2+}-\lambda_{2-})} & 0 & 1 \\
            0 & 0 & 1 & 0 \\
            0 & 0 & 0 & e^{n(\lambda_{4+}-\lambda_{4-})}
        \end{pmatrix}
        \quad \mbox{on $(-z_1,-z_2)$},\label{jumpT3}\\
        T_+ & = T_-
        \begin{pmatrix}
            1 & 0 & 0 & e^{n(\lambda_{4+}-\lambda_{1-})} \\
            0 & 1 & 0 & e^{n(\lambda_{4+}-\lambda_{2})} \\
            0 & 0 & 1 & 0 \\
            0 & 0 & 0 & 1
        \end{pmatrix}
        \quad \mbox{on $(-z_2,0)$}, \label{jumpT4} \\
        T_+ & = T_-
        \begin{pmatrix}
            1 & 0 & e^{n(\lambda_{3+}-\lambda_{1-})} & 0 \\
            0 & 1 & e^{n(\lambda_{3+}-\lambda_{2})} & 0 \\
            0 & 0 & 1 & 0 \\
            0 & 0 & 0 & 1
        \end{pmatrix}
        \quad \mbox{on $(0,z_2)$}, \label{jumpT5} \\
        T_+ & = T_-
        \begin{pmatrix}
            e^{n(\lambda_{1+}-\lambda_{1-})} & 0 & 1 & 0 \\
            0 & 1 & e^{n(\lambda_{3+}-\lambda_{2})} & 0 \\
            0 & 0 & e^{n(\lambda_{3+}-\lambda_{3-})} & 0 \\
            0 & 0 & 0 & 1
        \end{pmatrix}
        \quad \mbox{on $(z_2,z_1)$},\label{jumpT6}
    \end{align}
    \begin{align}
        T_+ & = T_-
        \begin{pmatrix}
            1 & 0 & e^{n(\lambda_3-\lambda_1)} & 0 \\
            0 & 1 & e^{n(\lambda_3-\lambda_2)} & 0 \\
            0 & 0 & 1 & 0 \\
            0 & 0 & 0 & 1
        \end{pmatrix}
        \quad \mbox{on $(z_1,x^*)$}, \label{jumpT7}\\
        T_+ & = T_-
        \begin{pmatrix}
            1 & 0 & e^{n(\lambda_3-\lambda_1)} &
            e^{n(\lambda_4-\lambda_1)} \\
            0 & 1 & e^{n(\lambda_3-\lambda_2)} &
            e^{n(\lambda_4-\lambda_2)} \\
            0 & 0 & 1 & 0 \\
            0 & 0 & 0 & 1
        \end{pmatrix}
        \quad \mbox{on $(x^*,+\infty)$}. \label{jumpT8}
    \end{align}
    On the vertical segment (with upwards orientation), we have
    \begin{align} \label{jumpT9}
        T_+ & = T_-
        \begin{pmatrix}
            1 & 0 & 0 & 0 \\
            0 & 1 & 0 & 0 \\
            0 & 0 & 0 & 1 \\
            0 & 0 & -1& e^{n(\lambda_{4+}-\lambda_{4-})}
        \end{pmatrix}
        \quad \mbox{on $(-iz_3,iz_3)$}.
    \end{align}
    The jumps on $\Sigma$ are as follows, where we take clockwise orientation on $\Sigma$, so
    that $T_+(z)$ ($T_-(z)$) for $z \in \Sigma$, is the limiting value of
    $T(z')$ as $z'\to z$ from outside (inside) of $\Sigma$:
    \begin{align} \label{jumpT10}
        T_+ & = T_-
        \begin{pmatrix}
            1 & 0 & 0 & 0 \\
            0 & 1 & 0 & 0 \\
            0 & 0 & 1 & 0 \\
            0 & 0 & e^{n(\lambda_3-\lambda_4)} & 1
        \end{pmatrix}
        \quad \mbox{on $\{z \in \Sigma \mid \Re z <0\}$},\\
        T_+ & = T_-
        \begin{pmatrix}
            1 & 0 & 0 & 0 \\
            0 & 1 & 0 & 0 \\
            0 & 0 & 1 & e^{n(\lambda_4-\lambda_3)} \\
            0 & 0 & 0 & 1
        \end{pmatrix}
        \quad \mbox{on $\{z \in \Sigma \mid \Re z >0\}$} \label{jumpT11}
    \end{align}
    \item[(3)] $T(z) = I+O(1/z)$ as $z \to \infty$.
\end{enumerate}

Now it follows that the jump matrix in (\ref{jumpT9})
on the vertical segment $(-iz_3,iz_3)$ tends to a constant matrix
since $\Re \lambda_{4+} < \Re \lambda_{4-}$. This may also be deduced
from Figures \ref{5figuur3a} and \ref{5figuur3b}, since we have
$\Re \lambda_{3} < \Re \lambda_4$ in the right half-plane within
the closed contour in these figures. The strict inequality remains
valid up to the vertical segment $(-iz_3,iz_3)$
so that $\Re \lambda_{3-} < \Re \lambda_{4-}$, since the
orientation of the segment is upwards. Since $\Re \lambda_{4 \pm} = \Re \lambda_{3\mp}$ the inequality
$\Re \lambda_{4+} < \Re \lambda_{4-}$ follows.

All the jump matrices on the real line
tend to the identity matrix as $n \to \infty$, except
for the ones on the intervals $(-z_1,-z_2)$ and $(z_2,z_1)$. The
1,1-entry and the 3,3-entry of the jump matrix on $(z_2,z_1)$,
see (\ref{jumpT6}), are rapidly oscillating for large $n$,
and the same holds for the 2,2-entry and
the 4,4-entry of the jump matrix on $(-z_1,-z_2)$, see
(\ref{jumpT3}). These oscillating entries are turned
into exponentially decaying ones in a standard
way by opening lenses around the intervals $[-z_1,-z_2]$ and
$[z_2,z_1]$. This will be done next.

\subsection{Third transformation: $T \mapsto S$}

\begin{figure}
    \begin{center}
        \includegraphics[height=2cm]{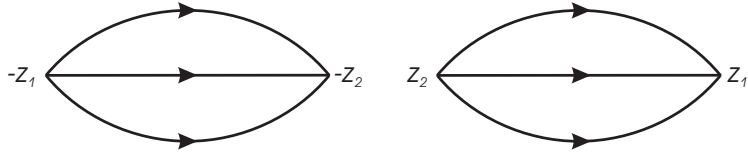}
        \caption{Opening of the lenses $\Gamma_1$ and $\Gamma_2$ around the intervals $[-z_1,-z_2]$
        and $[z_2,z_1]$.}\label{5figuur6}
    \end{center}
\end{figure}

We now turn the oscillatory entries on the diagonal of
(\ref{jumpT3}) and (\ref{jumpT6}) into
exponentially decaying ones. This is  done in a standard way
by opening lenses $\Gamma_1$ and $\Gamma_2$ in  (small) neighborhoods of
$(-z_1,-z_2)$ and $(z_2,z_1)$, respectively, as shown in Figure~\ref{5figuur6}.
Define $S$ as follows:
\begin{align}
    S & = T \quad \mbox{ outside $\Gamma_1$ and $\Gamma_2$},\\
    S & = T
    \begin{pmatrix}
        1 & 0 & 0 & 0 \\
        0 & 1 & 0 & 0 \\
        0 & 0 & 1 & 0 \\
        0 & -e^{n(\lambda_2-\lambda_4)} & 0 & 1
    \end{pmatrix}
    \quad \mbox{in the upper lens region around $[-z_1,-z_2]$},
    \end{align}
    \begin{align}
    S & = T
    \begin{pmatrix}
        1 & 0 & 0 & 0 \\
        0 & 1 & 0 & 0 \\
        0 & 0 & 1 & 0 \\
        0 & e^{n(\lambda_2-\lambda_4)} & 0 & 1
    \end{pmatrix}
    \quad  \mbox{in the lower lens region around $[-z_1,-z_2]$}, \\
    S & = T
    \begin{pmatrix}
        1 & 0 & 0 & 0 \\
        0 & 1 & 0 & 0 \\
        -e^{n(\lambda_1-\lambda_3)} & 0 & 1 & 0 \\
        0 & 0 & 0 & 1
    \end{pmatrix}
    \quad  \mbox{in the upper lens region around $[z_2,z_1]$},
    \label{defTuprrechts}\\
    S & = T
    \begin{pmatrix}
        1 & 0 & 0 & 0 \\
        0 & 1 & 0 & 0 \\
        e^{n(\lambda_1-\lambda_3)} & 0 & 1 & 0 \\
        0 & 0 & 0 & 1
    \end{pmatrix}
    \quad  \mbox{in the lower lens region around $[z_2,z_1]$}.
\end{align}

\begin{figure}
    \begin{center}
        \includegraphics[width=10cm,height=6cm]{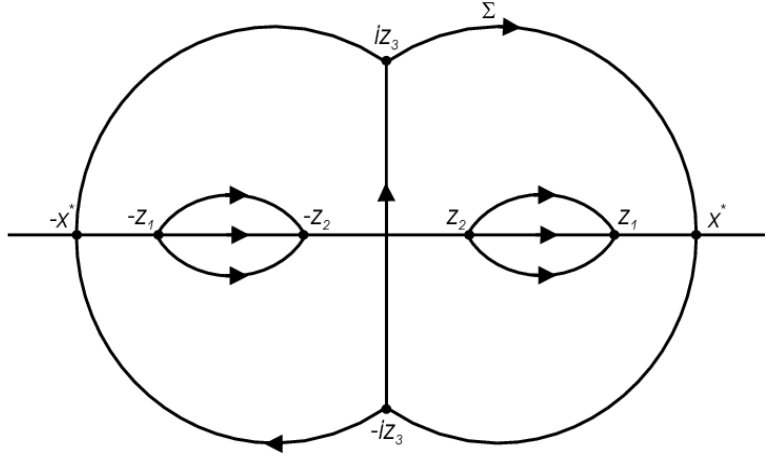}
        \caption{\label{5figuur3}
        Jump contour for the RH problem for $S$. The contour $\Sigma$ has clockwise
        orientation and the upper and lower lips of the lenses are oriented
        from left to right. The $4\times 4$ matrix valued
        function $S$ is analytic outside these contours, and has jumps (\ref{jumpS1})--(\ref{jumpS13})
        along the various parts of the contour.}
    \end{center}
\end{figure}
Then $S$ satisfies the following RH problem on the contour shown in Figure~\ref{5figuur3}:

\subsubsection*{RH problem for $S$}
\begin{enumerate}
    \item[(1)] $S$ is analytic on $\mathbb{C} \setminus (\mathbb{R}
    \cup \Sigma \cup \Gamma_1 \cup \Gamma_2 \cup [-iz_3,iz_3])$,
    \item[(2)] On the real line, $S$ has the following jumps:
    \begin{align}
        S_+ & = S_-
        \begin{pmatrix}
            1 & 0 & e^{n(\lambda_{3+}-\lambda_{1-})} &
            e^{n(\lambda_{4+}-\lambda_{1-})} \\
            0 & 1 & e^{n(\lambda_{3+}-\lambda_{2-})} &
            e^{n(\lambda_{4+}-\lambda_{2-})} \\
            0 & 0 & 1 & 0 \\
            0 & 0 & 0 & 1
        \end{pmatrix}
         \quad \mbox{on $(-\infty,-x^*)$}, \label{jumpS1} \\
        S_+ & = S_-
        \begin{pmatrix}
            1 & 0 & 0 & e^{n(\lambda_{4+}-\lambda_{1-})} \\
            0 & 1 & 0 & e^{n(\lambda_{4+}-\lambda_{2-})} \\
            0 & 0 & 1 & 0 \\
            0 & 0 & 0 & 1
        \end{pmatrix}
        \mbox{on $(-x^*,-z_1)$}, \label{jumpS2} \\
        S_+ & = S_-
        \begin{pmatrix}
            1 & -e^{n(\lambda_{2+}-\lambda_{1-})} & 0 & e^{n(\lambda_{4+}-\lambda_{1-})} \\
            0 & 0 & 0 & 1 \\
            0 & 0 & 1 & 0 \\
            0 & -1 & 0 & 0
        \end{pmatrix}
        \quad \mbox{on $(-z_1,-z_2)$},\label{jumpS3}\\
        S_+ & = S_-
        \begin{pmatrix}
            1 & 0 & 0 & e^{n(\lambda_{4+}-\lambda_{1-})} \\
            0 & 1 & 0 & e^{n(\lambda_{4+}-\lambda_{2})} \\
            0 & 0 & 1 & 0 \\
            0 & 0 & 0 & 1
        \end{pmatrix}
        \quad \mbox{on $(-z_2,0)$}, \label{jumpS4}
        \end{align}
        \begin{align}
        S_+ & = S_-
        \begin{pmatrix}
            1 & 0 & e^{n(\lambda_{3+}-\lambda_{1-})} & 0 \\
            0 & 1 & e^{n(\lambda_{3+}-\lambda_{2})} & 0 \\
            0 & 0 & 1 & 0 \\
            0 & 0 & 0 & 1
        \end{pmatrix}
        \quad \mbox{on $(0,z_2)$}, \label{jumpS5}\\
        S_+ & = S_-
        \begin{pmatrix}
            0 & 0 & 1 & 0 \\
            -e^{n(\lambda_{1+}-\lambda_{2})} & 1 & e^{n(\lambda_{3+}-\lambda_{2})} & 0 \\
            -1 & 0 & 0 & 0 \\
            0 & 0 & 0 & 1
        \end{pmatrix}
        \quad \mbox{on $(z_2,z_1)$},\label{jumpS6} \\
        S_+ & = S_-
        \begin{pmatrix}
            1 & 0 & e^{n(\lambda_{3}-\lambda_1)} & 0  \\
            0 & 1 & e^{n(\lambda_3-\lambda_2)} & 0 \\
            0 & 0 & 1 & 0 \\
            0 & 0 & 0 & 1
        \end{pmatrix}
        \quad \mbox{on $(z_1,x^*)$}, \label{jumpS7} \\
        S_+ & = S_-
        \begin{pmatrix}
            1 & 0 & e^{n(\lambda_3-\lambda_1)} &
            e^{n(\lambda_4-\lambda_1)} \\
            0 & 1 & e^{n(\lambda_3-\lambda_2)} &
            e^{n(\lambda_4-\lambda_2)} \\
            0 & 0 & 1 & 0 \\
            0 & 0 & 0 & 1
        \end{pmatrix}
        \quad \mbox{on $(x^*,+\infty)$}. \label{jumpS8}
    \end{align}
    On the vertical segment, we have that
    \begin{align}
        S_+ & = S_-
        \begin{pmatrix}
            1 & 0 & 0 & 0 \\
            0 & 1 & 0 & 0 \\
            0 & 0 & 0 & 1 \\
            0 & 0 & -1 & e^{n(\lambda_{4+}-\lambda_{4-})}
        \end{pmatrix}
        \quad \mbox{on $(-iz_3,iz_3)$.} \label{jumpS9}
    \end{align}
    On $\Sigma$, the jumps of $S$ are
    \begin{align}
        S_+ & = S_-
        \begin{pmatrix}
            1 & 0 & 0 & 0 \\
            0 & 1 & 0 & 0 \\
            0 & 0 & 1 & 0 \\
            0 & 0 & e^{n(\lambda_3-\lambda_4)} & 1
        \end{pmatrix}
        \quad \mbox{on $\Sigma$ in the left half-plane}, \label{jumpS10} \\
        S_+ & = S_-
        \begin{pmatrix}
            1 & 0 & 0 & 0 \\
            0 & 1 & 0 & 0 \\
            0 & 0 & 1 & e^{n(\lambda_4-\lambda_3)} \\
            0 & 0 & 0 & 1
        \end{pmatrix}
        \quad \mbox{on $\Sigma$ in the right half-plane}. \label{jumpS11}
    \end{align}
    The jumps on the lenses $\Gamma_1$ and $\Gamma_2$ around $[-z_1,-z_2]$ and
    $[z_2,z_1]$ are
    \begin{align}
        S_+ & = S_-
        \begin{pmatrix}
            1 & 0 & 0 & 0 \\
            0 & 1 & 0 & 0 \\
            0 & 0 & 1 & 0 \\
            0 & e^{n(\lambda_2-\lambda_4)} & 0 & 1
        \end{pmatrix}
        \quad \mbox{on $\Gamma_1$}, \label{jumpS12} \\
        S_+ & = S_-
        \begin{pmatrix}
            1 & 0 & 0 & 0 \\
            0 & 1 & 0 & 0 \\
            e^{n(\lambda_1-\lambda_3)} & 0 & 1 & 0\\
            0 & 0 & 0 & 1
        \end{pmatrix} \label{jumpS13}
        \quad \mbox{on $\Gamma_2$}.
    \end{align}
    \item[(3)] $S = I+O(1/z)$ as $z \to \infty$.
\end{enumerate}

Now it can be checked that all non-constant entries in the
jump matrices for $S$ tend to $0$ as $n \to \infty$. This could
be done by completing the Figures \ref{5figuur3a} and \ref{5figuur3b} by drawing
all curves where $\Re \lambda_j = \Re \lambda_k$ for $j \neq k$
from which one can deduce the sign of $\Re (\lambda_j - \lambda_k)$
in various regions in the complex plane.

For the jumps on $(-z_1,-z_2)$ and $(z_2,z_1)$ one can also
use a Cauchy-Riemann type argument. For the interval $(z_2,z_1)$
this is based on the fact that
$(\lambda_1 - \lambda_3)_+$ is purely imaginary on $(z_2,z_1)$
and its derivative is $(\xi_1-\xi_3)_+ = 2 i \Im \xi_{1+}$
where $\Im \xi_{1+} > 0$. Then by the Cauchy-Riemann equations we have
that  $\Re(\lambda_1-\lambda_3) < 0$  in a neighborhood above $(z_2,z_1)$.

\subsection{Parametrix away from the branch points}
We are now
going to solve the model RH problem, where we ignore all
exponentially small entries in the jump matrices for $S$.
So we only keep the jumps on $(-z_1,-z_2)$, $(z_2,z_1)$,
and $(-iz_3,iz_3)$ (which are also the cuts for the Riemann
surface) and we look for a $4 \times 4$ matrix
valued function $N$ that satisfies the following RH
problem:

\subsubsection*{RH problem for $N$}
\begin{enumerate}
    \item[(1)] $N$ is analytic on $\mathbb{C}\setminus ([-z_1,-z_2] \cup
    [z_2,z_1] \cup [-iz_3,iz_3])$,
    \item[(2)] $N$ satisfies the following jumps along the cuts:
    \begin{align}
        N_+ & = N_-
        \begin{pmatrix}
            1 & 0 & 0 & 0 \\
            0 & 0 & 0 & 1 \\
            0 & 0 & 1 & 0 \\
            0 & -1 & 0 & 0
        \end{pmatrix}
        \quad \mbox{on $(-z_1,-z_2)$}, \label{jumpN1}  \\
        N_+ & = N_-
        \begin{pmatrix}
            0 & 0 & 1 & 0 \\
            0 & 1 & 0 & 0 \\
            -1 & 0 & 0 & 0 \\
            0 & 0 & 0 & 1
        \end{pmatrix}
        \quad \mbox{on $(z_2,z_1)$}, \label{jumpN2} \\
        N_+ & = N_-
        \begin{pmatrix}
            1 & 0 & 0 & 0 \\
            0 & 1 & 0 & 0 \\
            0 & 0 & 0 & 1 \\
            0 & 0 & -1 & 0
        \end{pmatrix}
        \quad \mbox{on $(-iz_3,iz_3)$}, \label{jumpN3}
    \end{align}
    \item[(3)] as $z \to \infty$, we have that
    \begin{align}
        N(z) \to I+O(1/z).
    \end{align}
\end{enumerate}
By means of the rational parametrization of the
Riemann surface we can transform this problem to the
complex $v$-plane and thereby solve this RH problem explicitly.

\begin{figure}
    \begin{center}
        \includegraphics[width=10cm,height=6cm]{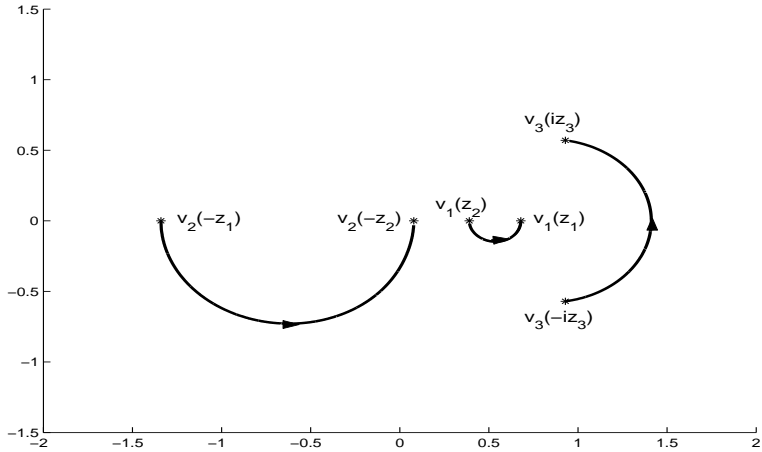}
        \caption{\label{5figuurN} Jump contours in the $v$-plane
        for the scalar RH problems for the functions $F_j$, $j=1,2,3,4$. }
    \end{center}
\end{figure}

Recall that we have the mapping $v = v_k(z)$ that maps the
$k$th sheet of the Riemann surface to the part $\omega_k$
of the $v$-plane, see also Figure~\ref{vplanet=0.25}.
We look for a solution for the RH problem for $N$ in
the form
\[ N(z) = \begin{pmatrix} F_j(v_k(z)) \end{pmatrix}_{j,k=1,\ldots,4} \]
where $F_j$, $j=1,\ldots, 4$ are four functions on the
$v$-plane.
Then the jump conditions on $N$ are satisfied if the $F_j$
are analytic except for the parts of the boundaries of the
domains $\omega_k$ that are indicated in Figure~\ref{5figuurN},
and on these parts we have $F_{j+} = - F_{j-}$.
The asymptotic condition on $N$ is satisfied provided that
\[ F_j(v_k(\infty)) = \delta_{jk}, \qquad \mbox{ for } j,k=1,\ldots, 4. \]
For each $j$, we then have a scalar RH problem for $F_j$ that
can be solved explicitly with elementary functions.

The precise form of $F_j$ is not important for what follows,
but we do need that $F_j(v)$ behaves like $(v-v_0)^{-1/2}$
as $v \to v_0$ and $v_0$ is one of the endpoints of the jump contour
in Figure~\ref{5figuurN}. For $N$ this implies that
\[ N(z) = O\left((z- z_0)^{-1/4}\right) \quad \mbox{ as } z \to z_0, \]
where $z_0$ is any of the branch points $\pm z_1$, $\pm z_2$, $\pm i z_3$.

\subsection{Parametrix near the branch points}
The jump matrices on $N$ and $S$ are not uniformly close to each
other near the branch points. That is why we need to treat these
points separately, and construct a local parametrix  $P$ around these
branch points.

We are going to construct a local parametrix around $z_1$. The local
parametrices around $-z_1$, $\pm z_2$, $\pm i z_3$ can be found in a similar
way, and are therefore not further discussed here. Consider a
small but fixed disc $U_{\delta}$ with radius $\delta$ around
$z_1$ that does not contain $z_2$. We then look for a $4 \times 4$
matrix valued function $P$ such that

\subsubsection*{RH problem for $P$ around $z_1$}
\begin{enumerate}
    \item[(1)] $P$ is analytic for  $z \in U_{\delta_0} \setminus (\Gamma_2 \cup
    \mathbb{R})$, for some $\delta_0> \delta$,
    \item[(2)] $P$ has the following jumps on the real line:
    \begin{align}
        P_+ & = P_-
        \begin{pmatrix}
            0 & 0 & 1 & 0 \\
            0 & 1 & 0 & 0 \\
            -1 & 0 & 0 & 0 \\
            0 & 0 & 0 & 1
        \end{pmatrix}
        \quad \mbox{on $U_{\delta} \cap [z_2,z_1]$}, \label{jumpP1} \\
        P_+ & = P_-
        \begin{pmatrix}
            1 & 0 & e^{n(\lambda_3-\lambda_1)} & 0 \\
            0 & 1 & 0 & 0 \\
            0 & 0 & 1 & 0 \\
            0 & 0 & 0 & 1
        \end{pmatrix}
        \quad \mbox{on $U_{\delta} \cap [z_1,+\infty)$}, \label{jumpP2}
    \end{align}
    and on the lens $\Gamma_2$ around $[z_2,z_1]$, we have
    \begin{align}
        P_+ & = P_-
        \begin{pmatrix}
            1 & 0 & 0 & 0 \\
            0 & 1 & 0 & 0 \\
            e^{n(\lambda_1-\lambda_3)} & 0 & 1 & 0 \\
            0 & 0 & 0 & 1
        \end{pmatrix}
        \quad \mbox{on $U_{\delta} \cap \Gamma_2$}. \label{jumpP3}
     \end{align}
    \item[(3)] As $n \to \infty$,
    \begin{align}
        P(z) = N(z)(I+O(1/n)) \quad \mbox{ uniformly for $z \in
        \partial U_{\delta} \setminus (\mathbb{R} \cup
        \Gamma_2)$}.\label{P1rand}
    \end{align}
\end{enumerate}
Note that the jumps for $P$ are not exactly the same as
the jumps for $S$. They differ by the entries
$e^{n(\lambda_{3+}-\lambda_2)}$ and $e^{n(\lambda_{1+}-\lambda_2)}$
which are exponentially small as $n \to \infty$.

The RH problem is solved in a standard way with Airy functions
\cite{DKMVZ1,DKMVZ2,ABK,BK2}.
 We have that
\begin{align}
    f_1(z) =
    \left[\frac{3}{4}(\lambda_1-\lambda_3)(z)\right]^{2/3}
\end{align}
is a conformal map that maps a neighborhood of $z_1$ onto a
neighborhood of the origin such that $f_1(z)$ is real and positive
for $z>z_1$.
We open the lens around $[z_2,z_1]$ such that $f_1$ maps
the part of $\Gamma_2$ in this neighborhood of $z_1$ to the rays with angles
$\frac{2\pi}{3}$ and $-\frac{2\pi}{3}$, respectively.
We put
\begin{equation} \label{defy012}
    y_0(s)=\Ai(s), \quad  y_1(s) = e^{2\pi i/3} \Ai(e^{2\pi i/3} s),
    \quad y_2(s) = e^{-2\pi i/3} \Ai(e^{-2\pi i/3} s),
\end{equation}
where $\Ai$ is the usual Airy function. Define the
matrix $\Psi$ by
\begin{align}
    \Psi(s) & =
    \begin{pmatrix}
        y_0(s) & 0 & -y_2(s) & 0 \\
        0 & 1 & 0 & 0 \\
        y_0'(s) & 0 & -y_2'(s) & 0 \\
        0 & 0 & 0& 1
    \end{pmatrix}
    \quad \mbox{for } \arg s \in (0, 2 \pi/3),\label{parz1region1} \\
    \Psi(s) & =
    \begin{pmatrix}
        -y_1(s) & 0 & -y_2(s) & 0 \\
        0 & 1 & 0 & 0 \\
        -y_1'(s) & 0 & -y_2'(s) & 0 \\
        0 & 0 & 0 & 1
    \end{pmatrix}
     \quad \mbox{for } \arg s \in (2\pi/3,\pi),\label{parz1region2} \\
    \Psi(s) & =
    \begin{pmatrix}
        -y_2(s) & 0 & y_1(s) & 0 \\
         0 & 1 & 0 & 0 \\
        -y_2'(s) & 0 & y_1'(s) & 0 \\
         0 & 0 & 0 & 1
    \end{pmatrix}
    \quad \mbox{for } \arg s \in (-\pi,-2\pi/3), \\
    \Psi(s) & =
    \begin{pmatrix}
        y_0(s) & 0 & y_1(s) & 0 \\
        0 & 1 & 0 & 0 \\
        y_0'(s) & 0 & y_1'(s) & 0 \\
        0 & 0 & 0 & 1
    \end{pmatrix}
    \quad \mbox{for } \arg s \in (-2\pi/3,0). \label{parz1region4}
\end{align}
Then, for any analytic prefactor $E$, we have that
\begin{align}
    P(z) = E(z)
    \Psi\left(n^{2/3}f_1(z)\right)\diag\left(e^{\frac{n}{2}(\lambda_1(z)-\lambda_3(z))},
    1, e^{-\frac{n}{2}(\lambda_1(z)-\lambda_3(z))},1\right)
\end{align}
satisfies the parts (1) and (2) of the RH problem for $P$.
If we choose $E$ as
\begin{align}
    E = \sqrt{\pi}N
    \begin{pmatrix}
        1 & 0 & -1 & 0 \\
        0 & 1 & 0 & 0 \\
        -i & 0 & -i & 0 \\
        0 & 0 & 0 & 1
    \end{pmatrix}
    \begin{pmatrix}
        n^{1/6}f_1^{1/4} & 0 & 0 & 0 \\
        0 & 1 & 0 & 0 \\
        0 & 0 & n^{-1/6}f_1^{-1/4} & 0 \\
        0 & 0 & 0 & 1
    \end{pmatrix}
\end{align}
then $E$ is analytic and the part (3) is satisfied as well.

\subsection{Fourth transformation: $S \mapsto R$}
\begin{figure}
    \begin{center}
        \includegraphics[height=5cm]{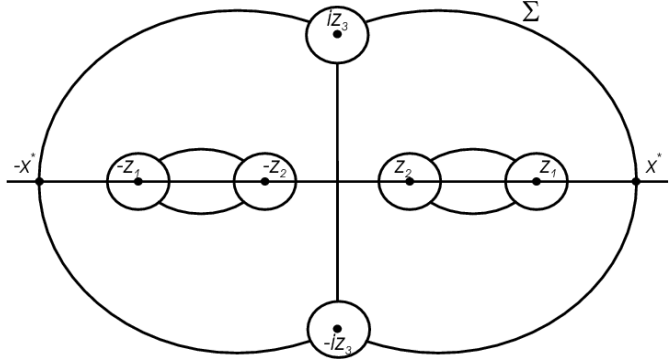}
        \caption{\label{figuur7}
        Jump contour for the RH problem for $R$.
        The $4\times 4$ matrix valued
        function $R$ is analytic outside these contours,
        and has jumps $R_+ = R_-(I + O(1/n))$
        uniformly on all parts of the contour.}
    \end{center}
\end{figure}

In the final transformation, we define the matrix valued function
$R$ as
\begin{align}
    R(z) & =  S(z)P(z)^{-1} \quad \mbox{in (small) discs around  $\pm z_1$, $\pm z_2$ and $\pm iz_3$}, \\
    R(z) & =  S(z)N(z)^{-1} \quad \mbox{outside the discs}.
\end{align}
Then $R$ is defined and is analytic (more precisely, has analytic continuation to the region)
outside the contour shown in Figure~\ref{figuur7}.

From the matching
condition (\ref{P1rand}) (and similar ones around the other
branch points), it follows that on
the circles around the branch points there is a jump
\begin{align}
    R_+ & = R_-(I+O(1/n)) \quad \mbox{ uniformly as $n \to
    \infty$}.
\end{align}
On the remaining contours, the jump is given by
\begin{align}
    R_+ & = R_-(I+O(e^{-cn})) \quad \mbox{ as $n \to
    \infty$}
\end{align}
for some constant $c>0$. Together with the asymptotic condition
\begin{align}
    R(z) = I+O(1/z) \quad \mbox{ as $z \to \infty$},
\end{align}
it then follows as in \cite{ABK,Dei,DKMVZ1,DKMVZ2,K1} that
\begin{align} \label{Rasymp}
    R(z) = I+O\left(\frac{1}{n(|z|+1)}\right) \quad \mbox{ as $n \to \infty$},
\end{align}
uniformly for $z \in \mathbb C \setminus \Gamma_R$, where $\Gamma_R$
is the jump contour for the RH problem for $R$, see Figure~\ref{figuur7}.

\section{Proofs of theorems for case $0 < t < t_{c,1}$}
\label{section5}

Having (\ref{Rasymp}) we can now prove Theorems \ref{theorem21}--\ref{theorem23}
in the same way as in \cite{BK2}.

\subsection{Proof of Theorem \ref{theorem21}}

Take $x,y \in (z_2,z_1)$. We follow the transformations $Y \mapsto U \mapsto T \mapsto S$,
to obtain from (\ref{defKinY}) that
\begin{align}
    K_n(x,y) =& \frac{e^{n(h(y)-h(x))}}{2 \pi i (x-y)}
    \begin{pmatrix}
        -e^{in\Im \lambda_{1+}(y)} & 0 & e^{-in\Im\lambda_{1+}(y)}
        & e^{n(\lambda_4(y) - \Re \lambda_{1+}(y))} \\
    \end{pmatrix} \nonumber
    \\ & \times S_+^{-1}(y) S_+(x)
    \begin{pmatrix}
        e^{-in \Im \lambda_{1+}(x)} \\ 0 \\ e^{in \Im
        \lambda_{1+}(x)} \\ 0
    \end{pmatrix}, \label{KinS}
\end{align}
where
\[ h(x) = \frac{1}{2} \left(\Re \lambda_{1+}(x) + \Re \lambda_{3+(x)}\right)
     -\frac{x^2}{2(1-t)} \mbox{ for } x > 0. \]
As in \cite[Section 9]{BK2} we obtain from (\ref{Rasymp}) that
\begin{align}
    S_+^{-1}(y)S_+(x) = I+O(x-y) \quad \mbox{ as $y \to x$}
\end{align}
uniformly in $n$, and therefore we get that
\begin{align}
    K_n(x,y)  = e^{n(h(y)-h(x))}
    \left(
    \frac{\sin(n\Im\left(\lambda_{1+}(x)-\lambda_{1+}(y))\right)}{\pi (x-y)}+O(1)
    \right),\label{kernxydicht}
\end{align}
where the $O(1)$ holds uniformly in $n$. Letting $y \to x$ we find
\begin{align}
    K_n(x,x) = \frac{n}{\pi} \Im \xi_{1+}(x) + O(1).
\end{align}
For $x \in (-z_2,-z_1)$ we get the similar relation but with
$\xi_{1+}(x)$ replaced by $\xi_{2+}(x)$.
Thus (\ref{rhox}) holds with $\rho$ defined by
\begin{equation} \label{defrho}
    \rho(x)  = \left\{ \begin{array}{ll}
     \frac{1}{\pi} \Im \xi_{1+}(x) & \mbox{ for } x > 0, \\[10pt]
     \frac{1}{\pi} \Im \xi_{2+}(x) & \mbox{ for } x < 0.
     \end{array} \right.
\end{equation}

The further statements in Theorem \ref{theorem21} are  now
easy consequences of the properties of $\xi_1$ and $\xi_2$.

\subsection{Proof of Theorem \ref{theorem22}}

Let $x_0 \in (z_2,z_1)$ and take
\[ x = x_0 + \frac{u}{n \rho(x_0)}, \qquad y = x_0 + \frac{v}{n \rho(x_0)}. \]
Then for $n$ large enough, we have $x,y \in (z_2,z_1)$, so that
(\ref{KinS}) holds.
Thus
\[ \frac{1}{n\rho(x_0)} \hat{K}_n(x,y)
    = \frac{\sin(n(\Im \lambda_{1+}(x)- \Im \lambda_{1+}(y)))}{\pi (u-v)}
    + O\left(\frac{1}{n}\right). \]
This leads to (\ref{sinekernellimit}) as in \cite[section 9.2]{BK2}.

The proof for $x_0 \in (-z_2,-z_1)$ is similar.

\subsection{Proof of Theorem \ref{theorem23}}

Take $x = z_1 + \frac{u}{(c_1n)^{2/3}}$ and $y = z_1 + \frac{v}{(c_1n)^{2/3}}$.
If $u, v < 0$, then we can follow the transformations
$Y \mapsto U \mapsto T \mapsto S \mapsto R$, to find from (\ref{defKinY})
\begin{align} \nonumber
    \frac{1}{(c_1n)^{2/3}}\hat{K}_n(x,y) =&
    \frac{1}{2\pi i(u-v)}
    \begin{pmatrix}
        -1& 0 & 1 & e^{n\Re(\lambda_{4+}(y)-\lambda_{1+}(y))}
    \end{pmatrix} \\
    & \times \nonumber
    \Psi^{-1}_+\left(n^{2/3}f_1(y)\right) E_n^{-1}(y)R_+^{-1}(y)R_+(x)E_n(x) \\
    & \times \Psi_+\left(n^{2/3}f_1(x)\right)
    \begin{pmatrix}
        1 \\0 \\1 \\0
    \end{pmatrix}.
\end{align}
The term $e^{n\Re(\lambda_{4+}(y)-\lambda_{1+}(y))}$ is exponentially
small, and does not contribute to the limit. Then we can use
the arguments of \cite[section 9.3]{BK2} to obtain (\ref{airykernellimit1}).
Similar arguments give (\ref{airykernellimit1}) in case $u > 0$
and/or $v > 0$.

Likewise we get (\ref{airykernellimit2}).

\section{Steepest descent analysis and proofs of theorems for
case $t_{c,1} < t < t_{c,2}$}\label{section6}

The steepest descent analysis is somewhat different for the
case $t_{c,1} < t < t_{c,2}$, due to the different sheet structure
of the Riemann surface, see Figure~\ref{figuur9f}. However the main lines
in the proof remain the same.
We only point out that now we have
two pairs of purely imaginary branch points, $\pm iz_2$ and $\pm i z_3$.
In the second transformation we need to open two global lenses in order to
remove the exponentially increasing entries in the jump matrices.
After the small opening of a lens around $(-z_1,z_1)$ we then construct
local parametrices at the branch points with Airy functions.

\begin{figure}
    \begin{center}
        \includegraphics[height=5cm]{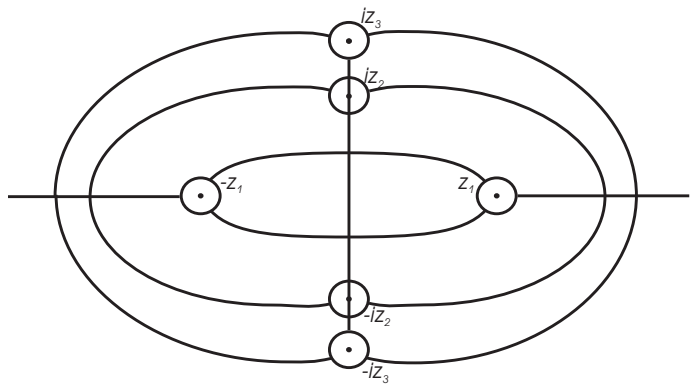}
        \caption{\label{figuur16}
        Jump contour for the RH problem for $R$ in case
        $t_{c,1} < t < t_{c,2}$ and $t \neq a/(a+b)$
        (for $t = a/(a+b)$ we have $z_2 = z_3$ and then
        we can simplify the contour by letting the two global
        lenses coincide).
        The $4\times 4$ matrix valued
        function $R$ is analytic outside these contours,
        and has jumps $R_+ = R_-(I + O(1/n))$
        uniformly on all parts of the contour.}
    \end{center}
\end{figure}

At the end of the transformations we arrive at a RH problem for $R$
with jumps on the contour shown in Figure~\ref{figuur16}. The jump
conditions for $R$ are  $R_+ = R_-(I + O(1/n))$, uniformly on all parts
of the contour. Then similar arguments and calculations lead to
the proofs of Theorems \ref{theorem21}--\ref{theorem23} for
the case $t_{c,1} < t < t_{c,2}$. We refer to \cite{Daems} for
complete details.

\end{document}